\documentclass[10pt]{iopart}

\usepackage{iopams}  

\expandafter\let\csname equation*\endcsname\relax

\expandafter\let\csname endequation*\endcsname\relax
\usepackage{pgfplots}
\usetikzlibrary{matrix,arrows,decorations.pathmorphing,patterns}
\pgfplotsset{compat=1.15}
\usepackage{mathrsfs}
\usetikzlibrary{arrows}
\usepackage{amsmath}
\usepackage{hyperref}
\newtheorem{theorem}{Theorem}[section] 

\newtheorem{corollary}[theorem]{Corollary}
\newtheorem{proposition}[theorem]{Proposition}

\newtheorem{remark}[theorem]{Remark}
\usepackage{tikz}
\definecolor{xdxdff}{rgb}{0.49019607843137253,0.49019607843137253,1.}
\definecolor{ududff}{rgb}{0.30196078431372547,0.30196078431372547,1.}
\definecolor{ffqqqq}{rgb}{1.,0.,0.}
\definecolor{uuuuuu}{rgb}{0.26666666666666666,0.26666666666666666,0.26666666666666666}
\definecolor{qqqqff}{rgb}{0.,0.,1.}
\begin{document}
	
	\title[]{Lipschitz Stability of an Inverse Problem of Transmission Waves with Variable Jumps}
	
	\author{L Baudouin$^1$, A Imba$^2$, 
		A Mercado$^2$ and A Osses$^3$}
	\address{$^1$CNRS ; LAAS ; 7 avenue du colonel Roche, F-31077 Toulouse, France and Université de Toulouse ; UPS, INSA, INP, ISAE, UT1, UTM, LAAS ; F-31077 Toulouse, France.}
	\address{$^2$ Departamento de Matem\'atica, Universidad T\'ecnica Federico Santa Mar\'{\i}a, Casilla 110-V, Valpara\'{\i}so, Chile}
	\address{$3$ Departamento de Ingenier\'ia Matem\'atica and Centro de Modelamiento Matem \'atico (UMI 2807 CNRS), FCFM Universidad de Chile, Casilla 170/3 - Correo 3, Santiago, Chile}
	
	\ead{lbaudoui@laas.fr, alex.imba@usm.cl, alberto.mercado@usm.cl, axosses@dim.uchile.cl}
	\vspace{10pt}
	\begin{indented}
		\item[]
	\end{indented}
	
	\begin{abstract}
		This article studies  an inverse problem for a transmission wave equation,
  a system where the main coefficient has a variable jump across an internal interface 
given by the boundary between two subdomains.
  The main result obtains   Lipschitz stability in recovering a zeroth-order coefficient
  in the equation. The proof is based on  the Bukhgeim-Klibanov method 
  and uses a new 
  one-parameter global 
  Carleman inequality,  specifically constructed  for the  case of a 
  variable main coefficient which is discontinuous  across a strictly convex interface.
   In particular, our hypothesis allows the main coefficient to vary smoothly within each subdomain up to the interface, thereby extending the preceding literature on the subject.
  	\end{abstract}
	
	%

	\section{Introduction}
	We consider an open and bounded set $\Omega \subset \mathbb{R}^d$, where $d \geq 2$, with a smooth boundary $\partial \Omega$, and a time $T > 0$. 
 We are interested in the
 wave equation given by
	\begin{equation}\label{wavevardis}
		\left\{   \begin{array}{rll}
			\partial^2_{t}u - \text{div}(a(x)\nabla u) + p(x)u = f  \quad &\text{in} \quad &(0, T) \times \Omega,\\
			u = g \quad &\text{on} \quad &(0, T) \times \partial \Omega, \\
			u(0,\cdot) = u^0, \quad \partial_t u(0,\cdot) = u^1 \quad &\text{in} \quad &\Omega,
		\end{array}
		\right. 
	\end{equation}
 where the \textit{potential} $p$ belongs to $L^\infty(\Omega)$ and the \textit{main coefficient} $a$ has the particularity of presenting a jump across an internal interface,
 in such a way that the system can be viewed as a {\it transmission wave equation}. 
 More precisely,   for a bounded set $\Omega_1 \subset \mathbb{R}^d$ satisfying $\overline{\Omega}_1 \subset \Omega$ and $\Omega_2 := \Omega \setminus \overline{\Omega}_1$, the main coefficient is given by 
	\begin{align}\label{variablecoef}
		a(x) =
		\begin{cases}
			a_1(x), &  x \in \Omega_1, \\
			a_2(x), & x \in \Omega_2,
		\end{cases}
	\end{align}
where each $a_1$ and $a_2$ is  smooth in $\Omega_1$ and $\Omega_2$, respectively. 

Therefore, the main  point of this framework is that the  main coefficient $a$ has a jump discontinuity on 
 $\Gamma_* := \partial \Omega_1$, 
 which  will be called the 
 \textit{interface}, 
  supposed to be strictly convex, in the sense that 
  the open set  $\Omega_1$ is a  strictly convex domain. We denote by $u_p$ the dependence of $u$ with respect $p$. The goal of this article is to investigate the stability in the following inverse problem: 
	\begin{quote}
 \textit{
 Given the boundary and source terms $g$ and $f$ and given the initial data $(u^0, u^1)$, can we determine the potential $p = p(x)$ from boundary observations of the flux of the solution $u_p$ through the boundary $\partial \Omega$?}
 \end{quote}


{ The process of determining $p=p(x)$ from flux measurements in a wave equation with a constant main coefficient has been extensively studied across the last decades. However, the same cannot be said for the case of discontinuous variable main coefficients, where variable jumps are present. This is precisely the main motivation behind this article. 
 
 Our approach for proving Lipschitz stability in this type of inverse problem is based on Carleman estimates and the Bukhgeim-Klibanov method, which was introduced in \cite{bukhgeim1981global}. This strategy has also been successfully applied to various inverse problems arising in wave equations that involve identifying sources \cite{puelGlobalEstimateLinear1996},\cite{puelGenericWellposednessMultidimensional1997}, potentials (zeroth order terms) \cite{yu.imanuvilovGlobalUniquenessStability2001}, main coefficients \cite{bellassoued*UniquenessStabilityDetermining2004}, \cite{imanuvilovDeterminationCoefficientAcoustic2003} \cite[Chapter 8]{isakovHyperbolicProblems2017} and other essential parameters within wave equations. The book by Bellassoued and Yamamoto \cite{bellassouedCarlemanEstimatesApplications2017} presents the applications of Carleman estimates for hyperbolic inverse problems in a systematic and clear manner. Additionally, the book by Klivanov and Timonov \cite{klibanovCarlemanEstimatesCoefficient2012} can be helpful for the numerical treatment of related inverse problems. One can also read \cite{baudouin2017convergent} for other approaches of the numerical reconstruction associated to wave equation inverse problems.}


 \subsection{Review of related   literature}  
	The study of inverse problems associated with the wave equation 
 with discontinuous main coefficient 
 is relevant to both real-world phenomena and from a theoretical perspective. 
From a mathematical viewpoint, it has attracted a lot of attention, 
in particular in control and inverse problems. 
Exact controllability issues for transmission waves with piecewise constant velocities were introduced by J. L. Lions in \cite{lionsHUM} using the duality between observability and controllability under typical geometric and time assumptions. Employing the multiplier method, Lions proved observability inequalities for the corresponding adjoint problem for inner domains that are star-shaped, assuming 
 that  the inner velocity is greater than the exterior one, 
 i.e., $a_1 \geq a_2$. The exact controllability results given by Lions can be easily adapted for waves with variable jumps under monotonicity of the traces 
 $$a_1(x)\geq a_2(x), \quad \forall x\in \Gamma_*,$$ and for star-shaped inner domains with respect to $x_0 \in \Omega_1$, with $a_1 \in W^{1,\infty}(\Omega_1)$, $a_2 \in W^{1,\infty}(\Omega_2)$, and the quite usual {\it multiplier condition}. Specifically, for $\rho \in (0,1]$, this condition is given by
	\begin{equation}\label{multiplierclassic}
		\nabla a \cdot (x - x_0) \leq 2a(1 - \rho) \quad \text{in}\quad \Omega_1 \cup \Omega_2. 
	\end{equation}
	
	It is well known that Carleman estimates for wave equations can 
 provide more powerful and more robust estimates than 
the  multiplier technique. 
 In this regard,  
in \cite{baudouinGlobalCarlemanEstimate2007}, it was obtained new  global two-parameter 
estimates and Lipschitz stability for potential recovery. More precisely, assuming the inner domain in $\mathbb{R}^2$  is strictly convex 
and the velocity possesses  a sign on the jump, 
a Carleman weight function  was constructed. 
The spatial part of the function defined in that work is essentially given by
 the square of the Minkowski functional of the inner subdomain,
 a function whose particular properties allow to obtain a convex function
with  non-vanishing gradient  and fullfiling  the same transmission conditions as the solution of the system. 
 Later, in \cite{baudouinInverseProblemSchrodinger2008}, 
 within the context  of a 
 Schr\"odinger equation, the construction  
 was extended to address the case of a transmission system in  $\mathbb{R}^d, d \geq 2$. 
It is noteworthy that the assumption of constant coefficients within each subdomain was essential for the approach presented in \cite{baudouinGlobalCarlemanEstimate2007}. This leaves open the question of whether the results can be extended to a broader class of discontinuous coefficients or geometries. In this work, we aim to address the case of variable jumps.

	Other works on Carleman estimates for transmission waves include: \cite{baudouinCarlemanEstimatesWave2022} for non-convex geometries of $\Omega_1$, \cite{filippasQuantitativeUniqueContinuation2022} proving quantitative unique continuation results (non Lipschitz-stability) employing a pseudodifferential approach, and \cite{gagnonSufficientConditionsControllability2023} using microlocal analysis. We mention also \cite{jiangStabilityInverseSource2023a} and \cite{zhaoInverseProblemTransmission2023a} for similar inverse problems. In \cite{riahiStabilityEstimateDetermination2015}, B. Riahi maintains the geometric and jump assumptions from \cite{baudouinGlobalCarlemanEstimate2007} and the study achieves H\"older stability in determining a discontinuous main coefficient. One of the conditions requires that $a_1$ and $a_2$ can vary on each subdomain, but they must remain constant in a neighborhood of the interface $\Gamma_*$.
	
	
	It has been observed that when trying to effectively reconstruct parameters, such as using the Carleman-based recontruction (Cb-Rec) algorithm presented in \cite{baudouinCarlemanBasedReconstructionAlgorithm2021}, it is advantageous to consider one-parameter Carleman estimates instead of two-parameter ones, mainly for numerical complexity reasons. Currently there are no references addressing inverse problems for transmission waves using one-parameter Carleman estimates. Therefore, this is a clear motivation to explore and establish such a novel Carleman estimates that can be effectively applied to ensure uniqueness, stability and reconstruction issues in inverse problems of waves with variable jumps. Theorem~\ref{StabThm} presents one first step towards addressing these questions.
	
	Our main objective in order to prove Theorem~\ref{StabThm} and to make completely precise the conditions under which it is valid, is to establish a new one-parameter global Carleman estimate for waves with variable jumps. Indeed, we aim to achieve this goal for a specific set of admissible coefficients as described in assumptions \eqref{infsup}-\eqref{rega} in Section~\ref{Bprop}. It is worth noting that this set of admissible coefficients goes beyond the results previously obtained by \cite{baudouinGlobalCarlemanEstimate2007}. Furthermore, we believe that our findings can also be applied to the results presented in \cite{riahiStabilityEstimateDetermination2015}. It is important to mention that we strive to provide a clearer understanding of the specific conditions under which our inverse problem and such a new Carleman estimate can be constructed for waves with variable jumps.
 
The rest of this introduction gives the construction of the main assumptions in Section~\ref{Bprop} and the main results of the article in Section~\ref{MR}.
	
	\subsection{Basic properties, weight function and admissible set of coefficients }\label{Bprop}
	
	Let us recall that if $a = a(x)$ given in \eqref{variablecoef} with $a_1 \in W^{1,\infty}(\Omega_1)$, $a_2 \in W^{1,\infty}(\Omega_2)$ is such that $$\inf_{x\in \Omega}a(x) \geq \alpha_0 > 0,$$ for some $\alpha_0 > 0$ and if $u^0 \in H^1(\Omega)$, $u^1 \in L^2(\Omega)$, $p \in L^\infty(\Omega)$, $g \in H^1(0,T;H^1(\partial \Omega))$, and $f \in L^1(0,T;L^2(\Omega))$, then, under the compatibility condition $g(0,\cdot) = u^0$ on $\partial \Omega $, we can prove that problem \eqref{wavevardis} admits a unique weak solution lying in the class
	\begin{equation*}
		u \in C([0,T]; H^1(\Omega)) \cap C^1([0,T]; L^2(\Omega)),
	\end{equation*}
	and  is such that $\partial_\nu u (t) \in L^2(\partial \Omega)$. This follows from \cite[Theorem 4.1]{lasiecka1986non}. 
 
 Moreover, assuming $\Omega_1$ is smooth, it is usual to consider an equivalent formulation of \eqref{wavevardis} in terms of $u_j = u_{\left|_{\Omega_j}\right.}$, $j=1,2$. The formulation in the corresponding subdomains is given by
	
	\begin{equation*}
		\left\{   \begin{array}{rll}
			\partial^2_{t}u_1 - \text{div}(a_1(x)\nabla u_1) + p(x)u_1 = f_1  \quad &\text{in} \quad &(0, T) \times \Omega_1,\\
			\partial^2_{t}u_2 - \text{div}(a_2(x)\nabla u_2) + p(x)u_2 = f_2  \quad &\text{in} \quad &(0, T) \times \Omega_2,\\
			u_2 = g \quad &\text{on} \quad &(0, T) \times \partial \Omega, \\
			u(0) = u^0(x), \quad \partial_t u(0) = u^1(x) \quad &\text{in} \quad &\Omega,
		\end{array}
		\right. 
	\end{equation*}
	with the transmission conditions:
	\begin{equation}
		u_1 = u_2 \quad {and} \quad a_1\frac{\partial u_1}{\partial \nu_1} + a_2\frac{\partial u_2}{\partial \nu_2} = 0 \quad \text{on} \quad (0,T) \times \Gamma_* \label{transuj}.
	\end{equation}
	Here, $\nu_1$ denotes the unit exterior normal vector to $\partial \Omega_1 = \Gamma_*$, and since $\partial \Omega_2 = \Gamma_* \cup \partial \Omega$, we can denote $\nu_2 = -\nu_1$. Additionally, $\nu$ denotes the normal exterior vector to $\partial \Omega$.
	
	In order to construct the new Carleman estimates of this work, we must introduce a suitable weight function and a set of assumptions on the main coefficients $a_1$ and $a_2$. To begin with, let us take $x_1 \in \Omega_1$ and consider the translated convex set $\Omega_1^{x_1} := \Omega_1 - x_1$. We introduce the Minkowski functional associated to $\Omega_1^{x_1}$, defined for each $x \in \mathbb{R}^d$ as follows:
	
	\begin{equation}\label{minkabs}
		p_{\Omega_1^{x_1}}(x) := \inf\left\lbrace h > 0: \frac{x}{h} \in \Omega_1^{x_1}\right\rbrace.
	\end{equation}
	Notice that $p_{\Omega_1^{x_1}}$ has several well-known properties which can be found for example in \cite[Proposition 5.1]{steinfunctional}. One of these properties is that it remains constant when restricted to the boundary of $\Omega_1^{x_1}$. We can make this property valid at the original $\Omega_1$ by composing $p_{\Omega_1^{x_1}}$ with $x-x_1$.  Moreover from \cite{baudouinGlobalCarlemanEstimate2007}, if $\Omega_1$ is strictly convex, then the square of Minkowski functional is strictly convex outside a ball of radius $\varepsilon>0$ small enough (to be specified later), that is, denoting $D^2$  the Hessian matrix:
 \begin{equation}\label{strictly:convex:function}
     \exists \  m_1>0 : D^2 \hspace{-0.05cm}\left[ p_{\Omega_1^{x_1}}(x -x_1) \right]^2(\xi, \xi) \geq m_1  |\xi|^2, \quad \forall \xi \in  \mathbb{R}^d\setminus\{0\}, \ \forall x \in (\Omega_1 \cup \Omega_2)\setminus B_{\varepsilon}(x_1).
 \end{equation}

Additionally,  given constants $0 < \alpha_0 \leq \alpha^0$ we suppose that the main coefficient $a = a(x)$, as given in equation \eqref{variablecoef}, fulfills the following
 \begin{quote}{\it - bounds conditions:}
	\begin{align}\tag{A.1}
		\alpha_0 \leq \inf_{x\in \Omega} a(x) \leq \sup_{x\in \Omega} a(x) \leq \alpha^0,\label{infsup} 
	\end{align}
   {\it- regularity conditions:}
	\begin{align}\tag{A.2}
		\begin{aligned}\label{reg} 
			&a_1 \in C^2(\overline{\Omega}_1) ,\\
			&a_2 \in C^2(\overline{\Omega}_2) .
		\end{aligned}
	\end{align}
  \end{quote}
	It has been noted that for obtaining Carleman estimates when $a_1$ and $a_2$ are constants, there is a sign restriction on the jump of the main coefficient. This jump condition is $a_1>a_2$ and the details can be found in reference \cite{baudouinGlobalCarlemanEstimate2007}. In our case, as $a_1$ and $a_2$ vary up to the interface, we need to introduce a sign assumption regarding the traces of $a_1$ and $a_2$ at the interface. The jump assumption we introduce writes as follows. 
 If we define the ratio 
 \begin{equation}\label{lambdadef}
		\lambda(x) := \frac{a_1(x)}{a_2(x)}, \quad x \in \Gamma_*,
	\end{equation} 
 we assume the following :
  \begin{quote}
	   {\it- jump condition:} {there exists a constant $\gamma>1$ such that }
\begin{equation}\label{jumpcond}
		\max_{x\in \Gamma*} \lambda(x) < \gamma < \min_{x\in \Gamma*}  \lambda^2(x). \tag{A.3}
	\end{equation}
  \end{quote}

 The assumption \eqref{jumpcond} arises precisely to obtain that certain integrals at the interface are positively dominant. More precisely, this assumption will arise for example for obtaining strictly positive terms of $s|w|^2$, $s|\partial_t w|^2$ and $s|\nabla_\tau w|^2$. We refer to the expressions in \eqref{J1} and \eqref{sumJ3J4J6} 
 concerning this point. 

Furthermore, in connection to assumption \eqref{multiplierclassic}, we introduce another assumption on $a_1$ and $a_2$ that typically arises when deriving global Carleman estimates for waves with a variable main coefficient. Given  $\left[ p_{\Omega_1^{x_1}}(x -x_1) \right]^2$, satisfying \eqref{strictly:convex:function} with $m_1$ and given  constants $\rho \in (0,1]$ and $\gamma$ fom \eqref{jumpcond}, we  set a
  \begin{quote}
    {\it- multiplier condition (with respect to $x_1 \in \Omega_1$ and $m_1>0$):}
	\begin{align}\tag{A.4}
		\begin{aligned}\label{rega}
		\nabla a (x) \cdot \nabla \left[ p_{\Omega_1^{x_1}}(x -x_1) \right]^2 \leq \frac{2}{\gamma}a(x)(1 - \rho)m_1, \quad x\in (\Omega_1 \cup \Omega_2) \setminus B_{\varepsilon}(x_1).
		\end{aligned}
	\end{align}
	 \end{quote}
  
	\begin{remark}
		It is worth noting that, 
  in the case that each one  of $a_1$ and $a_2$ is constant on the interface, we fix $\gamma=\frac{a_1}{a_2}$ and
then hypothesis   
  \eqref{jumpcond} is given  by $ \lambda =  \gamma <   \lambda^2$,
  it is equivalent to have $a_2 < a_1$. 
  This is the case considered in \cite{baudouinGlobalCarlemanEstimate2007}, where $a_1$ and $a_2$ are
  constants, and also of  \cite{riahiStabilityEstimateDetermination2015}, 
  where it is supposed that the coefficients are constants in a neighborhood of the interface.
  Therefore, assumption \eqref{jumpcond}  provides a sufficient condition for coefficients $a_1$ and $a_2$ that vary up to the interface, significantly expanding the range of coefficients for which Carleman estimates can be obtained.

	\end{remark}
	
	We define the weight function of our Carleman estimate as 
	\begin{equation}\label{varphiweight}
		\varphi(t,x)= \mu^{x_1}(x)- \beta t^2, \quad (t,x)\in \mathbb{R} \times \Omega,
	\end{equation} 
	which is a common approach (usually with $\mu^{x_1}(x) = |x-x_1|^2)$. We will explain how we construct the spatial part $\mu^{x_1}$ of our weight using a method similar to the one in \cite{baudouinGlobalCarlemanEstimate2007}, but adapted for variable jumps satisfying \eqref{jumpcond}. This construction will be detailed in the following paragraphs.
	
	  First, in view of  the singularity of the Minkowski functional at $x_1$, where the function is not regular, we introduce a  convenient cut-off function. 
 More precisely, for $\varepsilon_1$ and  $\varepsilon$ small enough (to be specified later) and such that  $0 < \varepsilon_1 < \varepsilon$, we consider a cut-off function 
 $\eta^{x_1} \in C^{\infty}(\Omega_1)$ such that
	\begin{eqnarray*}
 0 \leq \eta^{x_1} & \leq 1, \quad  x \in \overline{\Omega}_1, \\
		\eta^{x_1}(x) &= 0, \quad  x \in B_{\varepsilon_1}(x_1),\\
		\eta^{x_1}(x) &= 1, \quad  x \in \overline{\Omega}_1 \setminus \overline{B}_{\varepsilon}(x_1).
	\end{eqnarray*} 
	
	Then, using $\gamma>1$ given in \eqref{jumpcond} we finally define the spatial part of the weight as 
	\begin{equation}\label{weightinterior}
		\mu^{x_1}(x)=\left\{ \begin{array}{ll}
			\eta^{x_1}(x) \left[ p_{\Omega_1^{x_1}}(x -x_1) \right]^2  +M_1 \quad &x \in \overline{\Omega}_1,\\
			\gamma \left[ p_{\Omega_1^{x_1}}(x -x_1) \right]^2  + M_2, \quad &x \in\overline{\Omega}_2\setminus \overline{\Omega}_1,
		\end{array} \right.
	\end{equation}
	for $p_{\Omega_1^{x_1}}$ defined in \eqref{minkabs} and  constants
	\begin{equation}\label{constantsM1M2x1gamma}
		M_1>\gamma -1 \quad {and} \quad M_2:=M_1 -\gamma +1.
	\end{equation}
	
	\begin{remark}\label{remarkbmoweight}
		We highlight the dependence on $x_1$ by introducing it as a superscript. The function $\mu^{x_1}$, defined in \eqref{weightinterior}, follows 
  a similar  approach than   \cite{baudouinInverseProblemSchrodinger2008} and \cite{baudouinGlobalCarlemanEstimate2007}. 
		In these works, $\left[ p_{\Omega_1^{x_1}}(x -x_1) \right]^2 $  is expressed as
		\begin{equation}\label{minkbmo}
			\left[ p_{\Omega_1^{x_1}}(x -x_1) \right]^2  = \frac{|x-x_1|^2}{|y(x)-x_1|^2},
		\end{equation}
		where, for each $x \in \mathbb{R}^d \setminus \{x_1\}$, $y(x)$ is the unique point such that $y(x) \in \partial \Omega_1 \cap \ell(x_1,x)$, and $\ell$ denotes the segment joining $x_1$ and $x$; also, as 
  $a_1$ and $a_2$ are constants on $\Gamma_*$, the constant $\gamma$ is simply
  fixed as  $\gamma = \frac{a_1}{a_2}$. 
	\end{remark}
	
	The next result is established in \cite[Lemma 7]{baudouinGlobalCarlemanEstimate2007} and \cite[Proposition 8]{baudouinInverseProblemSchrodinger2008}. 
   
	\begin{proposition}\label{propositionweight}
		If $\Omega_1\subset \mathbb{R}^d$ is a $C^3$ strictly convex domain, then:
		\begin{eqnarray}
			\fl \exists \ \delta_1>0: \mu^{x_1}(x)\geq \delta_1 >0 \quad \forall x \in (\Omega_1 \cup \Omega_2)\setminus B_{\varepsilon}(x_1), \label{minkpositive}\\
			\fl \exists \  m_1>0 : D^2\hspace{-0.05cm}\,{\mu^{x_1}(x)}(\xi, \xi) \geq m_1  |\xi|^2, \quad \forall \xi \in  \mathbb{R}^d\setminus\{0\}, \quad \forall x \in (\Omega_1 \cup \Omega_2)\setminus B_{\varepsilon}(x_1), \label{minkhessian}\\
			\fl \exists\  c_1>0 : |\nabla \mu^{x_1}(x)|^2 \geq  c_1^2>0, \quad \forall x \in  (\Omega_1 \cup \Omega_2)\setminus B_{\varepsilon}(x_1). \quad  \label{minkgradnonzero}
		\end{eqnarray}
	\end{proposition}

 It is important to note that the function $\mu^{x_1}$ constructed here  satisfies, in addition to the properties  
 stated in  Proposition \ref{propositionweight},  convenient transmission conditions. Those conditions are local properties near the interface $\Gamma_*$ and will be discussed in detail in Proposition \ref{propositionweightinterfce}.

	\subsection{Main results}\label{MR}

  The following theorem gives the main result of this article. It answers both the uniqueness and the Lipschitz stability of the stated inverse problem.

	\begin{theorem}[Lipschitz Stability]\label{StabThm}
		Let $\overline{\Omega}_1 \subset \Omega$ be a $C^3$ strictly convex domain and $x_1,x_2 \in \Omega_1$ with $x_1\neq x_2$. Consider  $\Omega_1^{x_i}:= \Omega_1 - x_i$, for $i=1,2$ and $p_{\Omega_1^{x_i}}$ the Minkowski functional associated to $\Omega_1^{x_i}$ defined according \eqref{minkabs}. Let $0<\alpha_0\leq \alpha^0$ , $\rho \in(0,1]$, $m>0$. Assume  that the discontinuous main coefficient $a=a(x)$ given by \eqref{variablecoef} satisfies the bounds and regularity conditions detailed in \eqref{infsup} and \eqref{reg}. Further, suppose that $\lambda$ defined in \eqref{lambdadef}
 satisfies the jump condition \eqref{jumpcond} and the multiplier condition \eqref{rega} with respect to both $p_{\Omega_1^{x_i}}$, $i=1,2$ . Consider $\mu^{x_1}$ and  $\mu^{x_2}$ given by \eqref{weightinterior} allowing to define  
 \begin{equation}\label{Gamma1}
		\Gamma^{+}_{\mu^{x_i}} := \left\lbrace x \in \partial \Omega: \partial_\nu \mu^{x_i} > 0 \right\rbrace, \qquad i=1,2.
	\end{equation}	 
  Defining
	$
		L^\infty_{\leq m}(\Omega) = \left\lbrace p \in L^\infty(\Omega) : \|p\|_{L^\infty(\Omega)} \leq m \right\rbrace, ~
	$
 assume that $p\in L^\infty_{\leq m}(\Omega)$
 and the initial datum $u^0 \in H^1(\Omega)$ is such that 
		\begin{equation*}
			|u^0(x)| \geq \delta > 0, \ a.e \ {in} \ \Omega
		\end{equation*}
and  that the solution $u_p$ of \eqref{wavevardis} satisfies the regularity assumption 
		\begin{equation*}
			u_p \in H^1(0,T;L^\infty(\Omega)).
		\end{equation*}
		Then there exists $T_0 > 0$ 
  such that for any $T > T_0$ 
  there exists $C > 0$ such that, for any $q \in L^\infty_{\leq m}(\Omega)$, we have
		\begin{equation*}
			\|p - q\|_{L^2(\Omega)} \leq C \|\partial_{\nu} \partial_t (u_p - u_q)\|_{L^2(0,T;\Gamma_0)},
		\end{equation*}
  where $\Gamma_0 := \Gamma^+_{\mu^{x_1}} \cup \Gamma^+_{\mu^{x_2}}$.
	\end{theorem}

	The next result concerns Carleman estimates for the operator $L_q = (\partial^2_{t} - \text{div} (a(x)\nabla) + q)$, where $a = a(x)$ is given  in \eqref{variablecoef} and satisfies assumptions \eqref{infsup}-\eqref{rega}. We will use the notation $v(\pm T)$ to denote the evaluation of functions at $t=T$ and $t=-T$.  The precise statement of the Carleman estimates is the following result.
	\begin{theorem}[Carleman estimates]\label{mainthm} Assume $\Omega_1$ is a $C^3$ strictly convex domain. Let $x_1 \in \Omega_1$ and $p_{\Omega_1^{x_1}}$ the Minkowski functional associated to $\Omega_1^{x_1}:= \Omega_1 -x_1$, satisfying \eqref{strictly:convex:function} with $m_1>0$ and any $\varepsilon>0$.  Let  $0<\alpha_0\leq \alpha^0$ , $\rho \in(0,1]$, $m\geq 0$ and $\gamma>1$. Assume that $a=a(x)$ satisfies \eqref{infsup}, \eqref{reg}, \eqref{jumpcond} and \eqref{rega}. That is $a=a(x)$ belongs to
\begin{align*}
			\mathcal{A}_\gamma:=	\left\{
			\begin{array}{l}
				\left.a\right|_{\Omega_j} \in C^2, j=1,2. \quad \alpha_0 \leq \inf_{x\in \Omega} a(x) \leq \sup_{x\in \Omega} a(x) \leq \alpha^0, \\
				\sqrt \gamma a_2(x)   < a_1(x) < 
  \gamma a_2(x), \, \forall x \in \Gamma_*,  \\ \nabla a (x) \cdot \nabla \left[ p_{\Omega_1^{x_1}}(x -x_1) \right]^2 \leq \frac{2}{\gamma}a(x)(1 - \rho)m_1, \quad x\in (\Omega_1 \cup \Omega_2) \setminus B_{\varepsilon}(x_1) 
			\end{array}
			\right\},
		\end{align*}
 and  let $\mu^{x_1}$ defined in \eqref{weightinterior} and $\Gamma^+_{\mu^{x_1}}$ defined in \eqref{Gamma1}. Consider $\varphi(t,x)=\mu^{x_1}(x)-\beta t^2$. Then, for  any $q \in L^\infty_{\leq m}(\Omega)$, there exist a parameter $\beta_1>0$,
 such that $\forall \beta \in (0, \beta_1)$ there exists $s_0>0$ and $C>0$ such that 

		\begin{eqnarray}
			\fl s\int_{-T}^{T}\int_{\Omega} e^{2s\varphi} (|\partial_t v|^2 + a|\nabla v|^2 +s^2|v|^2)  \,dxdt \nonumber \\ \leq   C \int_{-T}^{T}\int_{\Omega} e^{2s\varphi}(L_q v)^2\,dxdt + Cs\int_{-T}^{T}\int_{\Gamma^+_{\mu^{x_1}}} e^{2s\varphi}\left|\partial_\nu v \right|^2 \,d\sigma dt \nonumber\\ \quad+~Cs\int\!\!\!\int_{\{(t,x) \in (-T,T)\times \Omega: \  \varphi<0\}}  e^{2s\varphi} (|\partial_t v|^2 +a|\nabla v|^2 + s^2|v|^2) \,dxdt \nonumber \\
			\quad+ ~Cs\int_{-T}^{T}\int_{B_\varepsilon(x_1)}e^{2s\varphi} (|\partial_t v|^2 +a|\nabla v|^2 + s^2|v|^2) \,dxdt \nonumber\\
			\quad + ~Cs\int_{\Omega}e^{2s\varphi(\pm T)}(|\partial_t v(\pm T)|^2 +|\nabla v(\pm T)|^2 +s^2|v(\pm T)|^2) \,dx \nonumber \\
			\quad +~Cs^3\int\!\!\!\int_{\{(t,x) \in (-T,T)\times\Gamma_*: \ \varphi<0\}} e^{2s\varphi}|v|^2\,d\sigma dt \nonumber \\
   \quad+~Cs^2\int_{\Gamma_*}e^{2s\varphi(\pm T)}|v(\pm T)|^2 d\sigma,\label{CarlBallR}
		\end{eqnarray}
		for all $ s \geq s_0$, for all $v  \in L^2(-T,T;H_0^1(\Omega))$ satisfying the transmission conditions \eqref{transuj} with $L_q v \in  L^2(-T,T;L^2(\Omega))$ and $\partial_\nu v \in L^2((-T,T)\times \partial \Omega))$.
	\end{theorem}
	
	
	
	
 
 We will apply the above Carleman estimates to obtain Lipschitz stability of our inverse problem.  
However, 
that is not a straightforward process due to some  additional observations on the right-hand side of  \eqref{CarlBallR},
In order to deal with this problem, we shall carefully combine these
estimates with  the inequalities provided by another weight function. 
More precisely, we define a function similar to \eqref{varphiweight} by
	\begin{equation}\label{phi}
		\phi(t, x) = \mu^{x_2}(x) - \beta t^2, \quad (t,x)\in \mathbb{R} \times \Omega,
	\end{equation}
	where $\mu^{x_2}$ is associated with some $x_2 \in \Omega_1$, with $x_2\neq x_1$ and is given by \eqref{weightinterior} and \eqref{constantsM1M2x1gamma}, (replacing $x_1$ by $x_2$ in the definition). Beware that Remark~\ref{remarkbmoweight} and Proposition~\ref{propositionweight} are thus also satisfied by $\mu^{x_2}$, and that \eqref{minkhessian} is satisfied with a constant $m_2$ (instead of $m_1$).\\
 
	
	Finally, the Carleman estimates with weights $\varphi$ and $\phi$ hold for $\beta \in (0, \beta_1)$ and $\beta \in (0,\beta_2)$, respectively. Denoting
	$\beta_0 = \min\{\beta_1, \beta_2\}$, and
	\begin{equation}\label{Ldef}
		L := \max\left\{\sup_{\Omega} \sqrt{\mu^{x_1}(x)}, \sup_{\Omega} \sqrt{\mu^{x_2}(x)}\right\},
	\end{equation}
 the assumption $T > T_0$  reads $T > L/\sqrt{\beta_0}$, where
  \begin{equation} \label{timecond}
	 T_0:= \frac{\max\left\{\sup_{\Omega} \sqrt{\mu^{x_1}(x)}, \sup_{\Omega} \sqrt{\mu^{x_2}(x)}\right\}}{\sqrt{\beta_0}}.
	\end{equation}
	\begin{remark}
		In \eqref{timecond}, estimating the required time for Lipschitz stability is generally not optimal. 
We can see that our estimate for $T_0$ depends on $L$ and $\beta_0$.
  This is an expected dependence, since $L$ is 
   given by the size of the domain and 
$\beta_0$ depends on $\alpha_0$, which is 
 the lower bound for the wave speed $a$.
However,  $\beta_0$  also depends on other  constants (see  \eqref{minimumbeta} below).
Providing a sharp estimate  of the minimal time in terms of size of the domain and the velocity $a$   would be a very interesting issue. 
	\end{remark}
	
	The next parts of the article are structured as follows: In Section \ref{ProofCarleman}, we present the proof of the main tool of this work, Theorem \ref{mainthm}. The proof  consists of four steps, 
 being the second one a crucial part of this work. 
 There, we use the assumption \eqref{jumpcond} to deal with interface integrals that contain traces of the coefficients $a_1=a_1(x)$ and $a_2=a_2(x)$, which vary up to $\Gamma_*$. In the same section, we will prove a Carleman estimate with  initial kinetic energy. Finally, in Section \ref{lipschitzsection}, we present the proof of Lipschitz stability for the potential recovery in transmission waves with variable jumps,
 which mainly focuses on two key steps: the construction of a convenient cut-off function and the application of the  Bukhgeim-Klibanov method.
	\section{Carleman estimates} \label{ProofCarleman}
	\subsection{A one-parameter global Carleman estimate for discontinuous wave operator}
	In this section, we aim to prove Theorem \ref{mainthm}. Our strategy is based on four main steps, following \cite{baudouinCarlemanEstimatesWave2022}. Firstly, we conduct computations in a generic domain while considering all the traces at the boundary resulting from the integration by parts. Secondly, we adapt the computations to the specific situation of variable jumps and carefully study the resulting terms at the interface. Thirdly, we establish a weighted norm in the interior of the domains, adopting the approach from \cite{baudouinCarlemanBasedReconstructionAlgorithm2021}. Finally, we combine the estimations from the previous steps and absorb negligible terms.
	\subsubsection*{Step 1 Generic  computations  }
	We denote the operator \(Lu = \partial_t^2u - \text{div}(a(x)\nabla u)\) and we perform the change of variables 
	\begin{eqnarray*}
		w &= e^{s\varphi}v, \quad s > 0,\\
		Pw &= e^{s\varphi} L(e^{-s\varphi}w).
	\end{eqnarray*}
	We will perform various computations in a generic set \(Q = (-T, T) \times U\) without making any assumptions regarding the boundary condition for \(w\) on \(\Sigma = (-T, T) \times \partial U\). At the same time, let us assume that $a$ is a smooth wave speed and that \(\mu^{x_1}\) in \(\varphi(t,x) = \mu^{x_1}(x) - \beta t^2\) is also smooth in such a way that all the upcoming computations make sense. Later on, we will see that it is sufficient to consider \(\mu^{x_1}\) given by \eqref{weightinterior}, and that the computations can be carried out even if $a$ is discontinuous. Computing \(Pw\), we write
	\begin{equation*}
		Pw = P_1w + P_2w + Rw,
	\end{equation*}
	where, for some function \(\alpha = \alpha(x)\) to be fixed later,
	\begin{eqnarray}
		P_1w &= \partial_t^2w - a\Delta w + s^2(|\partial_t \varphi|^2 - a|\nabla \varphi|^2)w - 2\nabla a \cdot \nabla w, \label{P1w}\\
		P_2w &= -2s\partial_t \varphi \partial_t w + 2sa\nabla \varphi \cdot \nabla w + \alpha as w, \nonumber\\
		Rw &= -s(\partial_t^2\varphi - \text{div}(a\nabla \varphi))w - \alpha as w + \nabla a \cdot \nabla w. \label{Rw}
	\end{eqnarray}
	As usual, we have
	\begin{align}\label{Pw2}
		\begin{aligned}
			\int_{-T}^{T}\int_{U}|P_1w|^2\,dxdt + \int_{-T}^{T}\int_{U}|P_2w|^2\,dxdt + 2(P_1w,P_2w)_{L^2(Q)} \\= \int_{-T}^{T}\int_{U} |Pw - Rw|^2\,dxdt.
		\end{aligned}
	\end{align}
	
	Let us compute $(P_1w,P_2w)_{L^2(Q)}$, and  denote by \(I_{i,j}\) the inner product of the \(i\)-th term of \(P_1w\) with the \(j\)-th term of \(P_2w\). Denoting 
 \(  \int\!\!\!\int_{Q} \,dxdt := \int_{-T}^{T} \int_{U} \,dxdt \) 
 and 
 \( \int\!\!\!\int_{\Sigma}\,d\sigma dt := \int_{-T}^{T} \int_{\partial U} \,d\sigma dt \),
 we obtain, using integration by parts in time, that
	\begin{eqnarray*}
		\fl I_{1,1} = -2s\int\!\!\!\int_{Q} \partial_t\varphi \partial_tw  \partial_t^2 w \,dxdt
		= s\int\!\!\!\int_{Q} |\partial_t w|^2 \partial_t^2\varphi \,dxdt - 
  \left. s\int_{U} |\partial_t w|^2 \partial_t \varphi \,dx\right\vert_{-T}^{T}.
	\end{eqnarray*}
	Similarly, integrating by parts in time and using that \(\partial_t \nabla \varphi = 0\), we get
	\begin{align*}
		&I_{1,2} = 2s\int\!\!\!\int_{Q} \partial_t^2w a \nabla \varphi \cdot \nabla w \,dxdt\\
		&= -2s\int\!\!\!\int_{Q} \partial_t w a\nabla \varphi \cdot \nabla \partial_t w \,dxdt + \left.2s\int_{\Omega}  \partial_twa\nabla \varphi \cdot \nabla w \right\vert_{-T}^{T}\,dx\\
		&= s\int\!\!\!\int_{Q}|\partial_t w|^2 \text{div}(a\nabla \varphi)\,dxdt - s\int\!\!\!\int_{\Sigma}  |\partial_t w|^2 a \frac{\partial \varphi}{\partial \nu} \,d\sigma dt + 
  \left. 2s\int_{U}  \partial_twa\nabla \varphi \cdot \nabla w \,dx\right\vert_{-T}^{T}.
	\end{align*}

	Noticing that \(a=a(x)\) and \(\alpha=\alpha(x)\) do not depend on time, we have
	\begin{align*}
		I_{1,3} = s\int\!\!\!\int_{Q}\alpha a\partial_t^2ww\,dxdt
		= -s \int\!\!\!\int_Q \alpha a|\partial_t w|^2\,dxdt + 
  \left. s \int_{U} \alpha a \partial_tww \,dx \right\vert_{-T}^{T},
	\end{align*}
	and in the same way, using integration by parts in space and time, we get
	\begin{align*}
		I_{2,1} &= 2s\int\!\!\!\int_{Q}a\Delta w \partial_t\varphi \partial_t w \,dxdt\\&=-2s\int\!\!\!\int_{Q}\nabla w \cdot \nabla (a\partial_t\varphi \partial_tw )\,dxdt + 2s\int\!\!\!\int_{\Sigma} \frac{\partial w}{\partial \nu}a\partial_t\varphi \partial_tw \,d\sigma dt \\
		&= s \int\!\!\!\int_{Q} a|\nabla w|^2 \partial_t^2\varphi \,dxdt - 2s \int\!\!\!\int_{Q} \partial_t\varphi \partial_tw \nabla a\cdot \nabla w \,dxdt - \left. s\int_{U}|\nabla w|^2 a \partial_t\varphi\,dx \right|_{-T}^{T}
		\\&\quad+ 2s \int\!\!\!\int_{\Sigma} a \frac{\partial w}{\partial \nu} \partial_t\varphi \partial_t w \,d\sigma dt.
	\end{align*}
	
	Integration by parts in space also gives
	\begin{align*}
		I_{2,2}&=-2s\int\!\!\!\int_{Q}a^2 \Delta w \nabla \varphi \cdot \nabla w \,dxdt\\
		&=2s \int\!\!\!\int_{Q}\nabla w \cdot \nabla (a^2 \nabla \varphi \cdot \nabla w) \,dxdt -2s\int\!\!\!\int_{\Sigma}\frac{\partial w}{\partial \nu}a^2 \nabla \varphi \cdot \nabla w \,d\sigma dt\\
		&=2s \int\!\!\!\int_{Q}a^2 D^2\, {\varphi}(\nabla w,\nabla w)\,dxdt +2s \int\!\!\!\int_{Q} \nabla w \cdot \nabla a^2 \nabla \varphi \cdot \nabla w \,dxdt\\ &\quad- s\int\!\!\!\int_Q |\nabla w|^2 \text{div}(a^2 \nabla \varphi) \,dxdt  -2s\int\!\!\!\int_{\Sigma} a^2 \nabla \varphi \cdot \nabla w \frac{\partial w}{\partial \nu}\,d\sigma dt \\ &\quad+s \int\!\!\!\int_{\Sigma} |\nabla w|^2 a^2 \frac{\partial \varphi}{\partial \nu} \,d\sigma dt,
	\end{align*}
	and 
	\begin{align*}
		I_{2,3}&=-s\int\!\!\!\int_{Q}\alpha a^2 \Delta w w \,dxdt=  s\int\!\!\!\int_{Q} \nabla w \cdot \nabla(\alpha a^2w)\,dxdt - s \int\!\!\!\int_{\Sigma}\alpha a^2 \frac{\partial w}{\partial \nu}  w \,d\sigma dt \\
		&=s\int\!\!\!\int_{Q} \alpha \nabla w \cdot \nabla a^2 w \,dxdt  +s \int\!\!\!\int_{Q}a^2 \nabla w\cdot \nabla \alpha w \,dxdt \\ &\quad + s\int\!\!\!\int_{Q}\alpha a^2 |\nabla w|^2 \,dxdt -s \int\!\!\!\int_{\Sigma}\alpha  a^2 \frac{\partial w}{\partial \nu} w \,d\sigma dt.
	\end{align*}
	
	Integration by parts in time, allows us to obtain
	\begin{align*}
		&I_{3,1}=-2s^3 \int\!\!\!\int_{Q}\partial_t \varphi \partial_t w (|\partial_t\varphi|^2 -a|\nabla \varphi|^2) w\,dxdt\\
		&= s^3 \int\!\!\!\int_{Q} |w|^2 \partial_t(\partial_t \varphi(|\partial_t\varphi|^2 -a|\nabla \varphi|^2)) \,dxdt 
  - \left. s^3\int_{U}|w|^2 \partial_t\varphi (|\partial_t \varphi|^2 -a|\nabla \varphi|^2)\,dx
  \right|_{-T}^{T} \\
		&= s^3 \int\!\!\!\int_{Q} |w|^2 [3\partial_t^2\varphi (|\partial_t \varphi|^2 -a|\nabla \varphi|^2) + 2a|\nabla \varphi|^2\partial_t^2\varphi]\,dxdt \\&\quad
  - \left. s^3\int_{U}|w|^2 \partial_t\varphi (|\partial_t \varphi|^2 
  -a|\nabla \varphi|^2)\,dx \right|_{-T}^{T},
	\end{align*}
	and similarly, integration in space gives
	\begin{align*}
		&I_{3,2}=2s^3 \int\!\!\!\int_{Q} a \nabla \varphi \cdot \nabla w (|\partial_t \varphi|^2 -a|\nabla \varphi|^2)w\,dxdt \\&= -s^3 \int\!\!\!\int_{Q} |w|^2 \text{div}(a\nabla \varphi (|\partial_t \varphi|^2 -a|\nabla \varphi|^2))\,dxdt \\&\quad +s^3 \int\!\!\!\int_{\Sigma} |w|^2a (|\partial_t\varphi|^2 -a|\nabla \varphi|^2)\frac{\partial \varphi}{\partial \nu}\,d\sigma dt\\
		&=-s^3 \int\!\!\!\int_{Q} |w|^2[\nabla a\cdot \nabla \varphi (|\partial_t\varphi|^2 -a|\nabla \varphi|^2)] \,dxdt \\&\quad-s^3 \int\!\!\!\int_{Q} |w|^2[a\Delta \varphi (|\partial_t\varphi|^2 -a|\nabla \varphi|^2) ]\,dxdt+2s^3\int\!\!\!\int_{Q}|w|^2a^2 D^2\,{\varphi}(\nabla \varphi, \nabla \varphi)\,dxdt \\&\quad+s^3 \int\!\!\!\int_{\Sigma} |w|^2a (|\partial_t\varphi|^2 -a|\nabla \varphi|^2)\frac{\partial \varphi}{\partial \nu}\,d\sigma dt.
	\end{align*}
	
	Direct computations give us
	\begin{equation*}
		I_{3,3}=s^3 \int\!\!\!\int_{Q}\alpha a \left(|\partial_t \varphi|^2 -a|\nabla \varphi|^2\right) |w|^2 \,dxdt,
	\end{equation*}
	and
	\begin{equation*}
		I_{4,1}= 4s \int\!\!\!\int_{Q}\partial_t\varphi \partial_tw \nabla a\cdot \nabla w \,dxdt. \\
	\end{equation*}
	
	Finally, note  also that
	\begin{equation*}
		I_{4,2}= -4s\int\!\!\!\int_{Q} a\nabla a\cdot \nabla w \nabla \varphi \cdot \nabla w \,dxdt= -2s\int\!\!\!\int_{Q}\nabla w\cdot \nabla a^2 \nabla \varphi \cdot \nabla w \,dxdt,
	\end{equation*}
	and
	\begin{equation*}
		I_{4,3}=-2s\int\!\!\!\int_{Q}\alpha a\nabla a\cdot \nabla w w \,dxdt =-s\int\!\!\!\int_{Q}\alpha \nabla w \cdot \nabla a^2 w \,dxdt. 
	\end{equation*}
	
	Gathering all the computed  terms and recalling $Q=(-T,T)\times U$, we write
	\begin{equation}\label{P1wP2wgeneric}
		(P_1w, P_2w)_{L^2(Q)}= A_{U} + Y_{U} + B_{\partial U}+E_{U}(\pm T),
	\end{equation}
	where $A_{U}$ is defined as  the sum of the so-called dominating interior terms, $B_{\partial U}$ is the sum of all the boundary terms , $Y_{U}$ contains a negligible interior term and $E_{U}(\pm T)$ in \eqref{evalinT} are integrals evaluated in $\pm T$. Taking into account the powers of $s$ we deduce
	\begin{align}\label{intintegrals}
		\begin{aligned}
			A_{U}=
			&s \int\!\!\!\int_{Q} |\partial_t w|^2 (\partial_t^2\varphi+ \nabla a\cdot \nabla \varphi +a\Delta \varphi -\alpha a )\,dxdt\\
			&+s \int\!\!\!\int_{Q} a|\nabla w|^2(\partial_t^2\varphi -2\nabla a\cdot \nabla \varphi - a\Delta \varphi + \alpha a)\,dxdt  \\
			&+2s \int\!\!\!\int_{Q}a^2 D^2\,{\varphi}(\nabla w,\nabla w)\,dxdt +2s \int\!\!\!\int_{Q} \partial_t \varphi \partial_t w\nabla a\cdot \nabla w \,dxdt   \\
			&+s^3 \int\!\!\!\int |w|^2 [(|\partial_ t\varphi|^2 -a|\nabla \varphi|^2)(3\varphi_{tt} -\nabla a\cdot \nabla \varphi - a\Delta\varphi +\alpha a)]\,dxdt \\&+ 2s^3 \int\!\!\!\int_{Q}|w|^2a^2D^2\,\varphi(\nabla \varphi, \nabla \varphi) \,dxdt +2s^3 \int\!\!\!\int_{Q}|w|^2a|\nabla \varphi|^2 \partial_t^2\varphi \,dxdt.
		\end{aligned}
	\end{align}
	The  negligible term is 
	\begin{align}\label{negligible}
		\begin{aligned}
			Y_{U}=s \int\!\!\!\int_{Q}a^2 \nabla w\cdot \nabla \alpha w\,dxdt,
		\end{aligned}
	\end{align}
	and the sum of the boundary terms is given by
	
	\begin{align}\label{boundaryterms}
		\begin{aligned}
			&B_{\partial U}\\
			&= -s\int\!\!\!\int_{\Sigma}  |\partial_t w|^2 a \frac{\partial \varphi}{\partial \nu} \,d\sigma dt + 2s \int\!\!\!\int_{\Sigma} a \frac{\partial w}{\partial \nu} \partial_ t\varphi \partial_ tw \,d\sigma dt -2s\int\!\!\!\int_{\Sigma} a^2 \nabla \varphi \cdot \nabla w \frac{\partial w}{\partial \nu} \,d\sigma dt\\ &\quad+s \int\!\!\!\int_{\Sigma} |\nabla w|^2 a^2 \frac{\partial \varphi}{\partial \nu} \,d\sigma dt- s \int\!\!\!\int_{\Sigma}\alpha a^2  \frac{\partial w}{\partial \nu}w\,d\sigma dt\\
			&\quad+s^3 \int\!\!\!\int_{\Sigma} |w|^2a (|\partial_t\varphi|^2 -a|\nabla \varphi|^2)\frac{\partial \varphi}{\partial \nu}\,d\sigma dt.
		\end{aligned}
	\end{align}
	
	Finally, the sum of the terms evaluated at $\pm T$ is
	
	\begin{align}\label{evalinT}
		\begin{aligned}
			&E_{U}(\pm T)\\ &=
   - \left. s\int_{U}|\partial_t w|^2 \partial_ t \varphi \,dx \right|_{-T}^{T} 
   +  \left. 2 s\int_{U} \partial_t wa\nabla \varphi \cdot \nabla w \,dx\right|_{-T}^{T} 
   + \left.  s \int_{U}\alpha a \partial_tw w \,dx \right|_{-T}^{T}\\
			&\quad - \left. s\int_{U}|\nabla w|^2 a \partial_t \varphi\,dx \right|_{-T}^{T} 
   - \left. s^3\int_{U}|w|^2 \partial_t \varphi (|\partial_t \varphi|^2 -a|\nabla \varphi|^2)\,dx \right|_{-T}^{T}.
		\end{aligned}
	\end{align}
	
	\subsubsection*{Step 2 At the interface}
	In this step, we first apply the computations performed in the previous step to each of the subdomains, considering $\Omega_1$ and $\Omega_2$ with the main coefficients given by $a_1=a_1(x)$ and $a_2=a_2(x)$ respectively. We then sum up the resulting terms. At this stage, we consider $\varphi(t,x) = \mu^{x_1}(x) - \beta t^2$ with $\mu^{x_1}$ given by \eqref{weightinterior}. From \eqref{P1wP2wgeneric}, we can express the inner product as follows:
	\begin{align*}
		\left(P_1 w, P_2 w\right)_{L^2((-T,T)\times \Omega)} &= A_{\Omega_1} + A_{\Omega_2} + Y_{\Omega_1} + Y_{\Omega_2} + B_{\partial \Omega_1} + B_{\partial \Omega_2} + E_{\Omega_1}(\pm T) \\ &\quad+ E_{\Omega_2}(\pm T).
	\end{align*}
	
	We write the sum of the boundary terms as follows
	\begin{align*}
		B_{\partial \Omega_1} + B_{\partial \Omega_2} = B_{\partial \Omega} + [B_{\Gamma_*}],
	\end{align*}
	where $[B_{\Gamma_*}]$ denotes the six integrals of \eqref{boundaryterms} supported on  $\Sigma_*= (-T,T)\times (\partial \Omega_1 \cap \partial \Omega_2)$, coming from the integration by parts in $\Omega_1$ and $\Omega_2$. More precisely, we introduce the notation
	\begin{align*}
		[B_{\Gamma_*}] = \sum_{i=1}^{6} (B_i(\Omega_1, \Gamma_*) + B_i(\Omega_2, \Gamma_*)) = \sum_{i
			=1}^{6}[J_i].
	\end{align*}
	
 We will bound by below $[B_{\Gamma_*}]$ using  the properties that $\mu^{x_1}$ satisfies at the interface. First, let us write $\mu^{x_1}$ given in \eqref{weightinterior} in the following way
	\begin{equation*}
		\mu^{x_1}(x) = \left\{ \begin{aligned}
			\mu^{x_1}_1(x), \quad &x \in \overline{\Omega}_1,\\
			\mu^{x_1}_2(x), \quad &x \in \overline{\Omega}_2\setminus \overline{\Omega}_1,
		\end{aligned} \right.
	\end{equation*}
	where $\mu^{x_1}_1=\left. \mu^{x_1}\right|_{\Omega_1}$ and $\mu^{x_1}_2= \left. \mu^{x_1}\right|_{\Omega_2}$. 
 The next result follows as a direct consequence of \cite[Proposition 2.1]{baudouinCarlemanEstimatesWave2022}.
	\begin{proposition}\label{propositionweightinterfce}
		If $\Omega_1$ is a $C^3$ strictly convex domain and $\gamma$ is defined in  \eqref{jumpcond}, then it holds $\mu_1^{x_1} \in C^{3}(\overline{\Omega}_1)$,  $\mu_2^{x_1} \in C^{3}(\overline{\Omega}_2)$ and
		\begin{align}
			\mu^{x_1}_1=\mu^{x_1}_2= cte, \quad  \text{on} \ \Gamma_*, \label{mucontGamma} \\
			\gamma \frac{\partial \mu^{x_1}_1}{\partial \nu_1} + \frac{\partial \mu^{x_1}_2}{\partial \nu_2}=0 , \quad  \text{on} \ \Gamma_*, \label{gammatransmission}\\
				\frac{\partial \mu^{x_1}_1}{\partial \nu_1} \geq\frac{2}{\text{diam}(\Omega)}, \quad   \text{on} \ \Gamma_*. \label{interfacemux1}
		\end{align}
	\end{proposition}
	\begin{remark}
		The transmission conditions \eqref{mucontGamma} and \eqref{gammatransmission} are proved in \cite{baudouinCarlemanEstimatesWave2022} and \cite{baudouinGlobalCarlemanEstimate2007} in the case of $\gamma=\frac{a_1}{a_2}$, with $a_1$ and $a_2$ constants at the interface. Here, we emphasize that even though $a_1$ and $a_2$ are not constants on $\Gamma_*$, we have  $\gamma>1$ as a constant thanks to \eqref{jumpcond}.
	\end{remark}
	Precisely, in $\overline{\Omega}_1\setminus \overline{B}_\varepsilon(x_1)$ we have
	\begin{align*}
		\mu^{x_1}_1(x)= \left[p_{\Omega_1^{x_1}}(x-x_1)\right]^2 +M_1
	\end{align*}
	and
	\begin{equation*}
		\Gamma_*= \left\{ x \in \mathbb{R}^d:  \quad p_{\Omega_1^{x_1}}(x-x_1)=1 \right\}. 
	\end{equation*}
	This last assertion comes from the definition of the Minkowski functional satisfying $p_{\Omega_1^{x_1}}(x-x_1)=1$ if and only if $x \in \Gamma_*$.
	In addition, in $\overline{\Omega}_2$ we have
	\begin{align*}
		\mu^{x_1}_2(x)= \gamma\left[p_{\Omega_1^{x_1}}(x-x_1)\right]^2 +M_2.
	\end{align*}
	Then, thanks to choosing $M_1$ and $M_2$ as in \eqref{constantsM1M2x1gamma} we readily obtain \eqref{mucontGamma} which reads
	\begin{align}\label{muctevalue}
		\mu^{x_1}_1(x)=1+M_1=\mu^{x_1}_2(x), \quad x\in \Gamma_*.
	\end{align}

	%
	%
	The transmission conditions \eqref{gammatransmission} and the condition \eqref{interfacemux1} follow by noting that $\Gamma_*$ is contained in the level curve $\left\{p_{\Omega_1^{x_1}}(x-x_1)=1\right\}$ and in this case we have
	\begin{align*}
		\nu_1= \frac{\nabla p_{\Omega_1^{x_1}}(x-x_1)}{\left|\nabla p_{\Omega_1^{x_1}}(x-x_1)\right|}=-\nu_2 \quad \text{and} \quad  \frac{\partial \mu^{x_1}_1}{\partial \nu_1}(x)= 2 \left|\nabla p_{\Omega_1^{x_1}}(x-x_1)\right|, \quad x\in \Gamma_*,
	\end{align*}
	with
	\begin{align*}
		\left|\nabla p_{\Omega_1^{x_1}}(x-x_1)\right|^2 \geq  \frac{\left[p_{\Omega_1^{x_1}}(x-x_1)\right]^2}{|x-x_1|^2} \geq \frac{1}{[\text{diam}(\Omega)]^2}, \quad x\in \Gamma_*.
	\end{align*}
	

	
	\begin{proposition}[Interface Integrals]\label{propinterface}
		Assume $\Omega_1$ is a $C^3$ strictly convex domain. Let $0 < \alpha_0 \leq \alpha^0$ such that $a = a(x)$ defined in \eqref{variablecoef} satisfies \eqref{infsup}, \eqref{reg} and \eqref{jumpcond}. Let $x_1 \in \Omega_1$, and let $\varphi(t,x) = \mu^{x_1}(x) - \beta t^2$ with $\mu^{x_1}$ satisfying \eqref{mucontGamma}, \eqref{gammatransmission}, and \eqref{interfacemux1}. If
		\begin{equation}\label{betaininterface}
			0<\beta < \frac{2 \alpha_0}{[\text{diam}(\Omega)]^2 (1+M_1)},
		\end{equation}
		then there exist $s_0 > 0$ and $C > 0$ such that 
		\begin{align}\label{minBGamma}
			-Cs^3\int\!\!\!\int_{\{(t,x) \in (-T,T)\times\Gamma_*: \ \varphi < 0\}} |w|^2\,d\sigma dt - Cs^2\int_{\Gamma_*} |w(\pm T)|^2 d\sigma \leq \sum_{j=1}^{6}[J_i] = [B_{\Gamma_*}],
		\end{align}
		for all $s \geq s_0$.
	\end{proposition}
	
	\begin{remark}
		Actually, it is possible to add, at the left-hand side of \eqref{minBGamma}, the terms
		\begin{align*}
			s\int_{-T}^{T}\int_{\Gamma_*} \left(|\partial_t w|^2 + s^2|w|^2  +  \left|\nabla_{\tau}w\right|^2+ \left|\frac{\partial w_1}{\partial \nu_1}\right|^2\right)\,d\sigma dt,
		\end{align*}
		where $\nabla_\tau$ denotes the tangential first-order operator of $w$ over $\Gamma_*$. This follows directly from the proof of Proposition \ref{propinterface}. See \eqref{miniJ1J6complete}.
		
	\end{remark}

	\subsubsection*{Proof of Proposition \eqref{propinterface}} Let us recall the notation 
	\begin{align*}
		\lambda(x) = \frac{a_1(x)}{a_2(x)}, \quad x \in \Gamma_*,
	\end{align*}
	and  $\Sigma_* = (-T,T) \times \Gamma_*$. 
	
	Recalling also the notation $v_1 = \left.v\right|_{\Omega_1}$, $v_2 = \left.v\right|_{\Omega_2}$ and the transmission conditions for $v$ as in \eqref{transuj}
	\begin{align*}
		v_1 = v_2, \quad \text{and} \quad  a_1\frac{\partial v_1}{\partial \nu_1} + a_2\frac{\partial v_2}{\partial \nu_2} = 0, \quad &\text{on} \quad (0,T) \times \Gamma_*,
	\end{align*}
	we can use the fact that $w = e^{s\varphi}v$ and the transmission conditions \eqref{mucontGamma} and \eqref{gammatransmission} satisfied by $\mu^{x_1}$ to obtain
	\begin{align}
		\begin{aligned}\label{contrasw}
			w_1 &= w_2, \quad \text{and} \quad    
			a_1\frac{\partial w_1}{\partial \nu_1} + a_2\frac{\partial w_2}{\partial \nu_2} &= sa_2w(\lambda-\gamma)\frac{\partial \mu^{x_1}_1}{\partial \nu_1}, \quad (t,x) \in (-T,T)\times\Gamma_*.
		\end{aligned}
	\end{align}

	To illustrate how we deal with each term of $[B_{\Gamma_*}]$, we provide details for $[J_1]$. 
 The remaining terms are treated in a similar way. First, from \eqref{boundaryterms}, we write 
	\begin{align*}
		[J_1] &= -s \int\!\!\!\int_{\Sigma_*}  |\partial_t{w_1}|^2 a_1 \frac{\partial \varphi_1}{\partial \nu_1}\,d\sigma dt - s \int\!\!\!\int_{\Sigma_*}  |\partial_t{w_2}|^2 a_2 \frac{\partial \varphi_2}{\partial \nu_2}\,d\sigma dt.
	\end{align*}
	
	Using \eqref{contrasw} and the transmission conditions \eqref{gammatransmission}, we have
	\begin{align}
		[J_1] &= -s \int\!\!\!\int_{\Sigma_*}  |\partial_t w|^2 \left(a_1 \frac{\partial \mu_1^{x_1}}{\partial \nu_1} +a_2 \frac{\partial \mu^{x_1}_2}{\partial \nu_2}\right)\,d\sigma dt \nonumber \\
		&= s\int\!\!\!\int_{\Sigma_*} |\partial_t w|^2 (a_2 \gamma - a_1) \frac{\partial \mu^{x_1}_1}{\partial \nu_1}\,d\sigma dt = s\int\!\!\!\int_{\Sigma_*} a_2 |\partial_t w|^2  \frac{\partial \mu^{x_1}_1}{\partial \nu_1} (\gamma- \lambda)\,d\sigma dt. \label{J1}
	\end{align}
	
	We can obtain, using \eqref{gammatransmission}, \eqref{contrasw}, 
 the decomposition $\nabla w= \left(\nabla_{\tau} w \, \vec{\tau}+ \frac{\partial w}{\partial \nu}\vec{\nu}\right)$, and 
  some  straightforward computations, that
	\begin{align*}
		&[J_3]+[J_4]=  s\int\!\!\!\int_{\Sigma_*}  a_2^2\left|\frac{\partial w_1}{\partial \nu_1}\right|^2 \frac{\partial \mu^{x_1}_1}{\partial \nu_1}\lambda^2(\gamma -1) \,d\sigma dt\\
		&\quad+s^3 \int\!\!\!\int_{\Sigma_*} a_2^2|w|^2 \left(\frac{\partial \mu^{x_1}_1}{\partial \nu_1}\right)^3\gamma (\lambda-\gamma)^2 \,d\sigma dt + s \int\!\!\!\int_{\Sigma_*}a_2^2\left|\nabla_{\tau} w \right|^2\frac{\partial \mu^{x_1}_1}{\partial \nu_1}(\lambda^2 - \gamma)\,d\sigma dt\\
		&\quad-2s^2  \int\!\!\!\int_{\Sigma_*} a_2^2 w \frac{\partial w_1}{\partial \nu_1} \left|\frac{\partial \mu^{x_1}_1}{\partial \nu_1}\right|^2 \lambda \gamma (\lambda - \gamma) \,d\sigma dt
	\end{align*}
	and 
	\begin{align*}
		&[J_6]\\ &= s^3 \int\!\!\!\int_{\Sigma_*} a_2|w|^2 \frac{\partial \mu^{x_1}_1}{\partial \nu_1} |\partial_t\varphi|^2 (\lambda -\gamma)\,d\sigma dt +s^3 \int\!\!\!\int_{\Sigma_*} a_2^2 |w|^2\left(\frac{\partial \mu^{x_1}_1}{\partial \nu_1}\right)^3(\gamma^3 - \lambda^2)\,d\sigma dt.
	\end{align*}
	These expressions imply
	\begin{align}\label{sumJ3J4J6}
		\begin{aligned}
			&[J_6]+[J_3]+[J_4] = 
			s^3 \int\!\!\!\int_{\Sigma_*} a_2|w|^2 \frac{\partial \mu^{x_1}_1}{\partial \nu_1} |\partial_t\varphi|^2 (\lambda -\gamma)\,d\sigma dt \\
			&+ s^3 \int\!\!\!\int_{\Sigma_*} a_2^2|w|^2\left(\frac{\partial \mu^{x_1}_1}{\partial \nu_1}\right)^3 \left( (\gamma^3 - \lambda^2 ) + \gamma (\lambda - \gamma)^2\right)\,d\sigma dt \\&+ s\int\!\!\!\int_{\Sigma_*}  a_2^2\left|\frac{\partial w_1}{\partial \nu_1}\right|^2 \frac{\partial \mu^{x_1}_1}{\partial \nu_1}\lambda^2(\gamma -1) \,d\sigma dt+ s \int\!\!\!\int_{\Sigma_*}a_2^2\left|\nabla_{\tau} w \right|^2\frac{\partial \mu^{x_1}_1}{\partial \nu_1}(\lambda^2 - \gamma)\,d\sigma dt\\
			&-2s^2  \int\!\!\!\int_{\Sigma_*} a_2^2 w \frac{\partial w_1}{\partial \nu_1} \left|\frac{\partial \mu^{x_1}_1}{\partial \nu_1}\right|^2 \lambda \gamma (\lambda - \gamma) \,d\sigma dt.
		\end{aligned}
	\end{align}
 
Moreover, using Young's inequality and noticing that Assumption \eqref{jumpcond} implies 
 $$\gamma > \max_{x\in \Gamma_*}\lambda(x) \geq \min_{x\in \Gamma_*}\lambda(x) > 1,$$
 we can estimate the last term on the right-hand side of \eqref{sumJ3J4J6} as follows:
	\begin{align}\label{crossterminterface}
		\begin{aligned}
			&\left| -2s^2  \int\!\!\!\int_{\Sigma_*} a_2^2 w \frac{\partial w_1}{\partial \nu_1} \left|\frac{\partial \mu^{x_1}_1}{\partial \nu_1}\right|^2 \lambda \gamma (\lambda - \gamma) \,d\sigma dt \right| \\
			&\leq 2s^2\int\!\!\!\int_{\Sigma_*} a_2^2 |w|\left| \frac{\partial w_1}{\partial \nu_1}  \right|\left|\frac{\partial \mu^{x_1}_1}{\partial \nu_1}\right|^2 \lambda \gamma (\gamma -\lambda)d\sigma  dt\\
			&\leq s^3 \int\!\!\!\int_{\Sigma_*} a_2^2|w|^2\left|\frac{\partial \mu^{x_1}_1}{\partial \nu_1}\right|^3 \gamma^2(\gamma-\lambda) \,d\sigma dt +s \int\!\!\!\int_{\Sigma_*} a_2^2\left|\frac{\partial w_1}{\partial \nu_1} \right|^2  \frac{\partial \mu^{x_1}_1}{\partial \nu_1} \lambda^2(\gamma -\lambda) d\sigma 
			dt\\
			&\leq s^3 \int\!\!\!\int_{\Sigma_*} a_2^2|w|^2\left|\frac{\partial \mu^{x_1}_1}{\partial \nu_1}\right|^3 \gamma^2(\gamma-1) \,d\sigma dt +s \int\!\!\!\int_{\Sigma_*} a_2^2\left|\frac{\partial w_1}{\partial \nu_1} \right|^2  \frac{\partial \mu^{x_1}_1}{\partial \nu_1} \lambda^2(\gamma -\lambda) d\sigma 
			dt.
		\end{aligned}
	\end{align}

Directly from  \eqref{crossterminterface}, \eqref{sumJ3J4J6},  \eqref{jumpcond} and \eqref{interfacemux1} we get that there exists a constant $C>0$ depending in particular on $\alpha_0$, $\min_{x \in \Gamma_*} \lambda(x)$, $\gamma$ and $\text{diam}(\Omega)$ and  changing from line to line such that
	\begin{align*}
		&[J_1]+ [J_3]+[J_4]+[J_6] \geq\\
		& Cs\int\!\!\!\int_{\Sigma_*} \left( |\partial_t w|^2   +  \left|\nabla_{\tau}w\right|^2 + \left|\frac{\partial w_1}{\partial \nu_1}\right|^2 \right)\,d\sigma dt\\
		&+s^3\int\!\!\!\int_{\Sigma_*}a_2^2 |w|^2(\gamma - \lambda)	\frac{\partial \mu^{x_1}_1}{\partial \nu_1} \left[\frac{4}{[\text{diam}(\Omega)]^2}(\gamma(\gamma - \lambda)+ \lambda + \gamma) - \frac{ |\partial_t \varphi|^2}{a_2} \right] \,d\sigma dt.
	\end{align*}
	Note also that since $|\partial_t \varphi|^2=4\beta^2 t^2$, we can deduce, using \eqref{muctevalue}, that 
	\begin{align*}
		&\Big\{(t,x) \in \Sigma_*: \varphi>0\}
		=\{(t,x) \in \Sigma_*: \mu^{x_1}(x)>\beta t^2\Big\} \\
		&=\Big\{(t,x) \in \Sigma_*: 4\beta \mu^{x_1}(x) > |\partial_t \varphi|^2 \Big\}
		=\Big\{(t,x) \in \Sigma_*: 4\beta (1+M_1) > |\partial_t \varphi|^2 \Big\}.
	\end{align*}
	In this last region, using \eqref{infsup}, \eqref{jumpcond} and \eqref{interfacemux1} we obtain 
	\begin{align*}
		&s^3\int\!\!\!\int_{\{(t,x) \in \Sigma_*: \varphi>0\}}a_2^2 |w|^2(\gamma - \lambda)	\frac{\partial \mu^{x_1}_1}{\partial \nu_1} \left[\frac{4}{[\text{diam}(\Omega)]^2}(\gamma(\gamma - \lambda)+ \lambda + \gamma) - \frac{ |\partial_t \varphi|^2}{a_2} \right] \,d\sigma dt\\
		&\geq s^3\int\!\!\!\int_{\{ \varphi>0\}} a_2^2 |w|^2 (\gamma-\lambda)\frac{\partial \mu^{x_1}_1}{\partial \nu_1} \left[ \frac{4}{[\text{diam}(\Omega)]^2}\left( \min_{x \in \Gamma_*} \lambda(x) + \gamma\right) -\frac{4\beta(1+M_1)}{\alpha_0} \right] \,d\sigma dt \\
		&\geq 4s^3\int\!\!\!\int_{\{\varphi>0\}} a_2^2 |w|^2 \left(\gamma-\lambda \right) \frac{\partial \mu^{x_1}_1}{\partial \nu_1} \left[ \frac{1}{[\text{diam}(\Omega)]^2}\left( 2\min_{x \in \Gamma_*} \lambda(x) \right) -\frac{\beta(1+M_1)}{\alpha_0} \right] \,d\sigma dt \\
		&\geq 4s^3\int\!\!\!\int_{\{(t,x) \in \Sigma_*: \varphi>0\}} a_2^2 |w|^2 \left(\gamma-\lambda \right) \frac{\partial \mu^{x_1}_1}{\partial \nu_1} \left[ \frac{2}{[\text{diam}(\Omega)]^2} -\frac{\beta(1+M_1)}{\alpha_0} \right] \,d\sigma dt. 
	\end{align*}
	Now, recalling the condition \eqref{betaininterface} on $\beta$ that reads
	$0<\beta<\frac{2 \alpha_0}{[\text{diam}(\Omega)]^2 (1+M_1)},$
	we obtain after using again \eqref{infsup}, \eqref{jumpcond}  and \eqref{interfacemux1} that
	\begin{multline}
	    \label{miniJ1J6complete}
	\relax[J_1]+ [J_3]+[J_4]+[J_6] \geq \,Cs\int\!\!\!\int_{\Sigma_*} \left( |\partial_t w|^2 +s^2|w|^2  +  \left|\nabla_{\tau}w\right|^2 +\left|\frac{\partial w_1}{\partial \nu_1}\right|^2 \right)\,d\sigma dt \\
   \quad - Cs^3\int\!\!\!\int_{\{(t,x) \in \Sigma_*: \varphi<0\}} |w|^2 \,d\sigma dt.
	\end{multline}
	
	Finally, it is not difficult to check that
	\begin{align*}
		\left| 	[J_2] +	[J_5] \right|  \leq Cs^2\int\!\!\!\int_{\Sigma_*} |w|^2 \,d\sigma dt + C\int\!\!\!\int_{\Sigma_*}\left|\frac{\partial w_1}{\partial \nu_1}\right|^2 \,d\sigma dt + Cs^2\int_{\Gamma_*} |w(\pm T)|^2 d\sigma.
	\end{align*}
	Then, there exists $s_0>0$ such that for all $s\geq s_0$
	\begin{multline}
	    \label{Prop23result}
	[B_{\Gamma_*}]= \sum_{j=1}^{6}[J_i]\geq Cs\int\!\!\!\int_{\Sigma_*} \left( |\partial_t w|^2 +s^2|w|^2  +  \left|\nabla_{\tau}w\right|^2+\left|\frac{\partial w_1}{\partial \nu_1}\right|^2\right)\,d\sigma dt\\ 
   \quad- Cs^3\int\!\!\!\int_{\{(t,x) \in \Sigma_*: \varphi<0\}} |w|^2 \,d\sigma dt -Cs^2\int_{\Gamma_*} |w(\pm T)|^2 d\sigma,
	\end{multline}
	and the proof of Proposition \ref{propinterface} is complete.

	\subsubsection*{Step 3 Interior strictly positive terms}
	
	Here, we work on the minimization of the interior integrals $A_{\Omega_1} + A_{\Omega_2}$. Our approach follows \cite{baudouinCarlemanBasedReconstructionAlgorithm2021}. 
 Taking into account that the interface is a zero-measure subset 
 of $\Omega$, each term gathered in $A_{\Omega_1} + A_{\Omega_2}$ 
 can be seen as an integral over $(-T,T)\times \Omega = Q$. 
 Also, we will refer to any derivative of $a$ as a function in $\Omega_1 \cup \Omega_2$, 
 using the notation $ \|a\|_{W^{1, \infty}} = \|a\|_{W^{1, \infty}(\Omega_1\cup\Omega_2)}$, and
 the same for other spaces.
 
 We denote by $A_j, \quad j=1, \ldots, 7$ the seven integrals in \eqref{intintegrals}. First, let us note that
	
	\begin{align*}
		A_1= \int\!\!\!\int_{Q} |\partial_t w|^2 \left( -2\beta + \nabla a \cdot \nabla \mu^{x_1} + a\Delta \mu^{x_1} - \alpha a \right) \,dxdt
	\end{align*}
	and
	\begin{align*}
		A_2+A_3= \int\!\!\!\int_{Q} a|\nabla w|^2 \left( -2\beta - 2\nabla a \cdot \nabla \mu^{x_1} - a\Delta \mu^{x_1} + \alpha a + 2a \frac{D^2\, \mu^{x_1}(\xi,\xi)}{|\xi|^2} \right)  \,dxdt,
	\end{align*}
	for all $\xi \in \mathbb{R}^d \setminus \left\lbrace 0 \right\rbrace $. 
 Notice that the function multiplying $|\nabla w|^2$ 
 in the above integral can be bounded by below 
  in  $(-T,T)\times Q\setminus B_{\varepsilon}(x_1)$ using \eqref{minkhessian} and \eqref{rega}. Indeed, as soon as  we set 
	$\alpha=\alpha(x)$ by
	\begin{equation}\label{fixalpha}
		\alpha = \Delta \mu^{x_1} - m_1 + \frac{3}{2a}\nabla a\cdot \nabla \mu^{x_1},
	\end{equation}
	and since $\mu^{x_1}$ and $a$ defined in \eqref{weightinterior} and \eqref{variablecoef}  are regular in each one of the sets $\overline{\Omega}_1$ and $\overline{\Omega}_2$, we obtain, using \eqref{infsup} that there exists $C>0$ (changing from line to line)  depending in particular on  $\|\mu^{x_1}\|_{C^2}, \|a\|_{W^{1,\infty}}$ and $\beta>0$ such that 
	
	\begin{align}\label{A2A3prev}
		A_2+A_3 \geq (-2\beta+ \alpha_0 m_1\rho )\int\!\!\!\int_{Q} a|\nabla w|^2 \,dxdt - Cs\int_{-T}^{T}\int_{B_\varepsilon(x_1)}a|\nabla w|^2\,dxdt.
	\end{align}
	
	Bringing the definition of $\alpha$ in $A_1$ and as a consequence of using \eqref{rega} and \eqref{infsup} we also obtain that there exists $C>0$  such that 
	
	\begin{align}\label{A1prev}
		A_1\geq (-2\beta +\alpha_0m_1 \rho) s\int\!\!\!\int_{Q}|\partial_t w|^2 \,dxdt - Cs\int_{-T}^{T}\int_{B_\varepsilon(x_1)}|\partial_t w|^2 \,dxdt.
	\end{align}
	
	%
	%
	
	Combining \eqref{A2A3prev} and \eqref{A1prev} we get
	
	\begin{align}\label{A1A2A3}
		\begin{aligned}
			A_{1} + A_2 + A_3 &\geq (-2\beta +\alpha_0m_1 \rho) s\int\!\!\!\int_{Q}(|\partial_t w|^2 + a|\nabla w|^2) \,dxdt \\&\quad- Cs\int_{-T}^{T}\int_{B_\varepsilon(x_1)}(|\partial_t w|^2 +a|\nabla w|^2) \,dxdt.
		\end{aligned}
	\end{align}

	Next, we  deal with the term $A_4$ given by
	\begin{align*}
		A_4&=2s \int\!\!\!\int_{Q} \partial_t \varphi \partial_t w\nabla a\cdot \nabla w\,dxdt\\&= 2s \int_{-T}^{T}\int_{\Omega_1} {\partial_t \varphi_1}{\partial_t w_1}\nabla a_1\cdot \nabla w_1 \,dxdt + 2s \int_{-T}^{T}\int_{\Omega_2} {\partial_t \varphi_2}{\partial_t w_2}\nabla a_2\cdot \nabla w_2 \,dxdt .
	\end{align*}
	Let us notice that, using  \eqref{infsup} and \eqref{reg}, there exists $C>0$ depending on $\beta, T$, $\|a\|_{W^{1,\infty}}$ and $\alpha_0$  such that
	\begin{align}\label{A4varphineg}
		\left| 2s \int\!\!\!\int_{\{ \varphi <0\}} \partial_t\varphi \partial_t w\nabla a\cdot \nabla w \,dxdt\right|\leq C s \int\!\!\!\int_{\{ \varphi <0\}} |\partial_t w|^2 + a|\nabla w|^2 \,dxdt,
	\end{align}
	where $\left\lbrace \varphi<0\right\rbrace :=\{ (t,x) \in (-T,T)\times \Omega: \ \varphi(t,x) <0\}$.
	On the contrary, in the set $\left\lbrace \varphi>0 \right\rbrace :=\{ (t,x) \in (-T,T)\times \Omega:\ \varphi(t,x) >0\}$ i.e.,
	\begin{align*}
		\Big\{(t,x) \in (-T,T)\times \Omega: \mu^{x_1}(x)>\beta t^2\Big\}
		=\left\{ 2\sqrt{\beta}\sqrt{\mu^{x_1}(x)}>|\partial_t \varphi (t,x)|\right\},
	\end{align*}
	we get, after using Young's inequality,  \eqref{infsup}, \eqref{reg} and the definition of $L$ in \eqref{Ldef}:	
	\begin{align}\label{A4varphipos}
		\begin{aligned}
			&\left| 2s \int\!\!\!\int\displaylimits_{\{ \varphi >0\}} \partial_t\varphi \partial_t w\nabla a\cdot \nabla w \,dxdt\right|\leq 2s \int\!\!\!\int_{\{ \varphi >0\}} |\partial_t\varphi|| \partial_t w||\nabla a||\nabla w| \,dxdt\\
			&\leq \frac{4s\sqrt{\beta}L \|\nabla a\|_{L^\infty} }{\sqrt{\alpha_0}}\int\!\!\!\int_{\{ \varphi >0\}} |\partial_t w|\sqrt{a}|\nabla w| \,dxdt\\
			&\leq \frac{2s \sqrt{\beta}L \|\nabla a\|_{L^\infty} }{\sqrt{\alpha_0}}\int\!\!\!\int_{\{ \varphi >0\}} \left( |\partial_t w|^2 + a|\nabla w|^2 \right) \,dxdt.
		\end{aligned}
	\end{align}
	Denoting 
	\begin{align}\label{rdefinition}
		r:= \frac{L \|\nabla a\|_{L^\infty} }{\sqrt{\alpha_0}},
	\end{align}
	we obtain, after combination of \eqref{A1A2A3}, \eqref{A4varphineg} and \eqref{A4varphipos} 
	\begin{align*}
		&A_1+A_2+A_3+A_4\geq \left(-2\beta +\alpha_0m_1 \rho -2\sqrt{\alpha_0}\sqrt{\beta}r\right)s\int\!\!\!\int_{Q}(|\partial_t w|^2 + a|\nabla w|^2) \,dxdt\\ 
		&-Cs\int_{-T}^{T}\int_{B_\varepsilon(x_1)}(|\partial_t w|^2 +a|\nabla w|^2) \,dxdt -Cs \int\!\!\!\int_{\{  \varphi<0\}}   (|\partial_t w|^2 +a|\nabla w|^2) \,dxdt.\\
	\end{align*}
	Assuming that $\beta$ satisfies
	\begin{align}\label{betafirstorder}
		0<\beta < \alpha_0\left( \sqrt{\left(\frac{r}{4}\right)^2+ \frac{m_1\rho}{2}} -\frac{r}{2}\right)^2,
	\end{align}
	so that $-2\beta +\alpha_0m_1 \rho -2\sqrt{\alpha_0}\sqrt{\beta}r \geq C >0$, then there exists $C>0$ such that
	\begin{align}\label{A1toA4}
		\begin{aligned}
			&A_1+A_2+A_3+A_4\geq C s\int\!\!\!\int_{Q}(|\partial_t w|^2 + a|\nabla w|^2) \,dxdt\\ 
			&-Cs\int_{-T}^{T}\int_{B_\varepsilon(x_1)}(|\partial_t w|^2 +a|\nabla w|^2) \,dxdt -Cs \int\!\!\!\int_{\{ \varphi<0\}}   (|\partial_t w|^2 +a|\nabla w|^2) \,dxdt.
		\end{aligned}
	\end{align}
	
	Let us deal with the integrals $A_5$, $A_6$ and $A_7$. We write  the sum these integrals in \eqref{intintegrals}, as follows
	\begin{align}\label{A5A6A7}
		A_5+A_6+A_7=	s^3 \int\!\!\!\int_{Q} |w|^3 c_w\,dxdt,
	\end{align}
	where $c_w=c_w(t,x)$ is given by
	
	\begin{align*}
		\begin{aligned}
			c_w&= (|\partial_t \varphi|^2 -a|\nabla \varphi|^2)(3\partial_t^2\varphi -\nabla a\cdot \nabla \varphi -a\Delta \varphi +\alpha a) + 2a^2D^2\,\varphi(\nabla \varphi, \nabla \varphi)\\ &\quad +2a|\nabla \varphi|^2 \partial_t^2\varphi.
		\end{aligned}
	\end{align*}
	On one hand, since $\varphi$ and $a$ are regular in each one of the sets $\overline{\Omega}_1$ and $\overline{\Omega}_2$, it is clear that there exists $C>0$ such that
	\begin{equation*}
		\sup_{(-T,T)\times (\Omega_1 \cup \Omega_2)} |c_w(t,x)| \leq C.
	\end{equation*}
	On the other hand, since  $\varphi(t,x)=\mu^{x_1}(x)-\beta t^2$, we can minimize  $c_w$ in $\Omega\setminus B_\varepsilon(x_1)$, using \eqref{minkhessian}, \eqref{fixalpha} and \eqref{rega} as follows:
	\begin{align*}
		&c_w= (4\beta^2 t^2 -a|\nabla \mu^{x_1}|^2)(-6\beta -\nabla a \cdot \nabla \mu^{x_1} - a\Delta \mu^{x_1} +\alpha a) +2a^2 D^2\, {\mu^{x_1}}(\nabla \mu^{x_1}, \nabla \mu^{x_1})\\ &\quad\quad-4a\beta|\nabla \mu^{x_1}|^2\\
		&\geq (4\beta^2 t^2 -a|\nabla \mu^{x_1}|^2) \left(-6\beta +\frac{1}{2} \nabla a \cdot \nabla \mu^{x_1} -am_1 \right) +2a^2 m_1|\nabla \mu^{x_1}|^2 -4a\beta |\nabla \mu^{x_1}|^2  \\
		&\geq -4\beta^2 t^2\left( 6\beta + \frac{|\nabla a||\nabla \mu^{x_1}|}{2} +am_1\right) + a |\nabla \mu^{x_1}|^2\left(6\beta  -\frac{\nabla a\cdot \nabla \mu^{x_1}}{2} +3am_1 -4\beta \right) \\
		&\geq -4\beta^2 t^2\left( 6\beta + \frac{|\nabla a||\nabla \mu^{x_1}|}{2} +am_1\right) + a |\nabla \mu^{x_1}|^2\left(2\beta + am_1 (2+\rho)  \right).
	\end{align*}
	
	Moreover, using \eqref{minkgradnonzero}, \eqref{infsup}, \eqref{reg} and \eqref{Ldef}  we can continue minimizing  $c_w$  in the region 
	\begin{align*}
		\lbrace(t,x) \in (-T,T)\times \Omega\setminus &B_\varepsilon(x_1)
		: \  \varphi(t,x) >0 \rbrace \\ &= \left\{(t,x) \in (-T,T)\times \Omega\setminus B_\varepsilon(x_1)
		: \ -4\beta^2 t^2 \geq -4\beta \mu^{x_1}(x) \right\}
	\end{align*}
	
	as follows:
	
	\begin{align*}
		c_w\geq& -4\beta \mu^{x_1}\left(  6\beta + \frac{|\nabla a||\nabla \mu^{x_1}|}{2} +am_1\right) + a |\nabla \mu^{x_1}|^2\left(2\beta + am_1 (2+\rho)  \right)\\
		=&2\beta \left( a|\nabla \mu^{x_1}|^2 - \mu^{x_1} |\nabla a||\nabla \mu^{x_1}|\right)  -24\beta^2 \mu^{x_1} +am_1\left( (2+\rho) a|\nabla \mu^{x_1}|^2 -4\beta \mu^{x_1}\right) \\
		\geq& 2\beta(\alpha_0 c_1^2 - L^2\|\nabla a\|_{L^\infty}\|\nabla \mu^{x_1} \|_{L^\infty}) -24\beta^2 L^2 +am_1\left( (2+\rho) \alpha_0c_1^2 -4\beta L^2 \right).
	\end{align*}
	In particular, recalling the notation for $r$ in \eqref{rdefinition} and if 
	\begin{align}\label{betaw2intermediate}
		\beta< \alpha_0 \frac{(2+\rho) c_1^2}{4L^2},
	\end{align}
	we get
	\begin{align}
		c_w\geq& 2\beta(\alpha_0 c_1^2 - L^2\|\nabla a\|_{L^\infty}\|\nabla \mu^{x_1} \|_{L^\infty}) -24\beta^2 L^2 +\alpha_0m_1\left( (2+\rho) \alpha_0c_1^2 -4\beta L^2 \right)\nonumber\\
		=&2\beta \alpha_0 \left( c_1^2  -L\|\nabla \mu^{x_1} \|_{L^\infty} r\right) -24\beta^2 L^2 + \alpha_0^2 m_1 (2+\rho)c_1^2 - 4\beta \alpha_0 L^2 m_1\nonumber\\
		\geq& -24\beta^2L^2 -2\beta \alpha_0\left|c_1^2 -L\|\nabla \mu^{x_1} \|_{L^\infty} r  - 2L^2 m_1\right| +   \alpha_0^2 m_1 (2+\rho)c_1^2.  \label{cwminimization}
	\end{align}
	Now, noticing that \eqref{cwminimization} corresponds to a second order polinomial in $\beta$, we assert that there exists a positive constant $C>0$ such that
	\begin{align*}
		-24\beta^2L^2 -2\beta \alpha_0\left|c_1^2 -L\|\nabla \mu^{x_1} \|_{L^\infty} r  - 2L^2 m_1\right| +   \alpha_0^2 m_1 (2+\rho)c_1^2 \geq C >0,
	\end{align*}
	for all $\beta$ satisfiying
	
	\begin{equation}\label{betaw2}
		0<\beta < \alpha_0 \rho_0,
	\end{equation}
	where 
	\begin{multline}
	  \label{rhodef}
			\rho_0:=\frac{1}{24 L^2} \left[\left|c_1^2 -L\|\nabla \mu^{x_1} \|_{L^\infty} r  - 2L^2 m_1\right|^2 +24L^2m_1(2+\rho)c_1^2\right]^{1/2}\\ \quad- \left|c_1^2 -L\|\nabla \mu^{x_1} \|_{L^\infty} r  - 2L^2 m_1\right|.  
	\end{multline}
	Therefore, \eqref{A5A6A7} becomes
	\begin{multline}
	    \label{A5toA7}
     A_5+A_6+A_7\geq Cs^3\int\!\!\!\int_{Q} |w|^2 \,dxdt 
	-Cs^3\int_{-T}^{T}\int_{B_\varepsilon(x_1)}|w|^2 \,dxdt \\
   \quad -Cs^3 \int\!\!\!\int_{\{  \varphi<0\}}   |w|^2  \,dxdt.
	\end{multline}
	
	\subsubsection*{Step 4 Conclusion}
	On one hand, taking into account that $B_{\partial \Omega_1}+ B_{\partial \Omega_2}= B_{\partial \Omega} + [B_{\Gamma_*}]$  we have, thanks to the Dirichlet homogeneous boundary condition of $v=e^{-s\varphi }w $,
	\begin{align*}
		\nabla w= \partial_\nu w \,\vec{\nu}, \ \text{on} \ \partial \Omega
	\end{align*}
	and then we  obtain 
	\begin{align}\label{boundarymin}
		B_{\partial \Omega} \geq -Cs\int_{-T}^{T} \int_{\Gamma^+_{\mu^{x_1}}} \left| \partial_\nu w \right|^2 \,d\sigma dt,
	\end{align}
	where $\Gamma^+_{\mu^{x_1}}$ is defined in \eqref{Gamma1}.
	
	On the other hand, let us take $\beta_1$ such that the conditions \eqref{betaininterface}, \eqref{betafirstorder}, \eqref{betaw2intermediate} and \eqref{betaw2} hold, that is,
	\begin{equation}
	    \label{minimumbeta}
		0<\beta_1< \alpha_0 \min\left\lbrace \frac{2}{[\text{diam}(\Omega)]^2 (1+M_1)} ,\left( \sqrt{\left(\frac{r}{4}\right)^2+ \frac{m_1\rho}{2}} -\frac{r}{2}\right)^2, \frac{(2+\rho) c_1^2}{4L^2}, \rho_0 \right\rbrace
	\end{equation}
	where $\rho_0$ is defined in \eqref{rhodef}. Then, we can combine  inequality  \eqref{Prop23result} with estimates \eqref{A1toA4}, \eqref{A5toA7} and \eqref{boundarymin} to conclude that for all  $\beta \in (0,\beta_1)$ there exists $s_0>0$ and $C>0$ such that, for all $s\geq s_0$:
	\begin{align}\label{summaryofestimates}
		\begin{aligned}
			&A_{\Omega_1} + A_{\Omega_2}+ B_{\partial \Omega_1} + B_{\partial \Omega_2} 
			\geq C s\int\!\!\!\int_{Q}(|\partial_t w|^2 + a|\nabla w|^2 +s^2|w|^2) \,dxdt\\
			&-Cs\int_{-T}^{T}\int_{B_\varepsilon(x_1)}(|\partial_t w|^2 +a|\nabla w|^2 +s^2|w|^2) \,dxdt -Cs\int_{-T}^{T} \int_{\Gamma^+_{\mu^{x_1}}} \left| \partial_\nu w\right|^2 \,d\sigma dt \\&-Cs \int\!\!\!\int_{\{(t,x) \in (-T,T)\times \Omega: \  \varphi<0\}}   (|\partial_t w|^2 +a|\nabla w|^2 +s^2|w|^2) \,dxdt \\
			&- Cs^3\int\!\!\!\int_{\{(t,x) \in \Sigma_*: \  \varphi<0\}} |w|^2 \,d\sigma dt  -Cs^2\int_{\Gamma_*} |w(\pm T)|^2 d\sigma.
		\end{aligned}
	\end{align}

	Recalling the definition of $Y_{U}$ in \eqref{negligible}, we also need to deal with 
	\begin{align*}
		Y_{\Omega_1} + Y_{\Omega_2} = s\int_{-T}^{T}\int_{\Omega_1}a_1^2 w_1 \nabla w_1 \cdot \nabla \alpha_1 \,dxdt + s\int_{-T}^{T}\int_{\Omega_2}a_2^2 w_2 \nabla w_2 \cdot \nabla \alpha_2 \,dxdt.
	\end{align*}
	Notice that from \eqref{fixalpha} we have $\alpha=\Delta \mu^{x_1} -m_1 +\frac{3}{2a}\nabla a\cdot \nabla \mu^{x_1}$ and
	\begin{align*}
		\sup_{x\in \Omega_1 \cup \Omega_2} |\nabla \alpha (x)| \leq C,
	\end{align*}
	with $C$ depending in particular on $\|a_j\|_{C^2}$, $\|\mu_j^{x_1}\|_{C^3}$ and $m_1$. Indeed, using \eqref{infsup} and Young's inequality we obtain
	
	\begin{align}\label{estnegligible}
		\left| Y_{\Omega_1} + Y_{\Omega_2} \right| \leq C \int\!\!\!\int_{Q} a|\nabla w|^2\,dxdt + Cs \int\!\!\!\int_{Q} |w|^2\,dxdt.
	\end{align}
	Similarly, from the definition of $E_{U}(\pm T)$ in \eqref{evalinT} we readily get
	\begin{align*}
		\left| E_{\Omega_1}(\pm T) + 	E_{\Omega_2}(\pm T)\right| \leq Cs\int_{\Omega}(|\partial_t w(\pm T)|^2 +|\nabla w(\pm T)|^2 +s^2|w(\pm T)|^2) \,dx.
	\end{align*}
	
	We also have from $Rw=-s(\partial_t^2\varphi - \text{div}(a\nabla \varphi))w - \alpha as w + \nabla a \cdot \nabla w$ given in \eqref{Rw} 
	\begin{align}\label{estresto}
		\begin{aligned}
			\int\!\!\!\int_{Q}|Rw|^2 \,dxdt &= \int\!\!\!\int_{Q} \left| -s(\partial_t^2\varphi - \text{div}(a\nabla \varphi))w - \alpha as w + \nabla a \cdot \nabla w \right|^2 \,dxdt  \\
			&\leq  Cs \int\!\!\!\int_{Q} |w|^2\,dxdt + 2\int\!\!\!\int_{Q} \left| \nabla a \cdot \nabla w \right|^2 \,dxdt\\
			&\leq Cs \int\!\!\!\int_{Q} |w|^2\,dxdt  + \frac{2\|\nabla a\|_{L^\infty} }{\alpha_0}\int\!\!\!\int_{Q} a\left| \nabla w \right|^2 \,dxdt\\
			&\leq Cs \int\!\!\!\int_{Q} |w|^2\,dxdt  + C\int\!\!\!\int_{Q} a\left| \nabla w \right|^2 \,dxdt.
		\end{aligned}
	\end{align}
	
	Finally, taking into account \eqref{Pw2} and since we have estimated by below the expression
	\begin{align*}
		\left(P_1 w, P_2 w\right)_{L^2((-T,T)\times \Omega)} &= A_{\Omega_1} + A_{\Omega_2} + Y_{\Omega_1} + Y_{\Omega_2} + B_{\partial \Omega_1} + B_{\partial \Omega_2} + E_{\Omega_1}(\pm T)\\ &\quad + E_{\Omega_2}(\pm T),
	\end{align*}
	we can use  $s$ to absorb \eqref{estnegligible} and \eqref{estresto} by the positive dominant terms on the right-hand side of \eqref{summaryofestimates}. More precisely, we obtain the  inequality 
	\begin{align}\label{Carlemaninw}
		\begin{aligned}
			&s\int\!\!\!\int_{Q}(|\partial_t w|^2 + a|\nabla w|^2 +s^2|w|^2) \,dxdt + s^3\int\!\!\!\int_{\Sigma_*} |w|^2\,d\sigma dt + 	\int\!\!\!\int_{Q}|P_1w|^2\,dxdt  \\ \leq C&	\int\!\!\!\int_{Q}|Pw|^2 \,dxdt 
			+ Cs\int_{-T}^{T} \int_{\Gamma^+_{\mu^{x_1}}} \left| \partial_\nu w\right|^2 \,d\sigma dt +Cs^3\int\!\!\!\int_{\{(t,x) \in \Sigma_*: \  \varphi<0\}} |w|^2 \,d\sigma dt\\
			+&Cs \int\!\!\!\int_{\{(t,x) \in (-T,T)\times \Omega: \  \varphi<0\}}   (|\partial_t w|^2 +a|\nabla w|^2 +s^2|w|^2) \,dxdt  \\ 
			+&Cs\int_{\Omega}(|\partial_t w(\pm T)|^2 +|\nabla w(\pm T)|^2 +s^2|w(\pm T)|^2) \,dx + 	Cs^2\int_{\Gamma_*} |w(\pm T)|^2 d\sigma\\
			+&Cs\int_{-T}^{T}\int_{B_\varepsilon(x_1)}(|\partial_t w|^2 +a|\nabla w|^2 +s^2|w|^2) \,dxdt,
		\end{aligned}
	\end{align}
	for $(s\geq s_0)$, where $Pw = e^{s\varphi} L(e^{-s\varphi}w)$ and $L=\partial_t^2 -\text{div}(a(x)\nabla)$.

	Returning to $v=e^{-s\varphi}w$ and using
	\begin{align*}
		|L v|^2 \leq C |L_q v|^2 + Cm|v|^2,
	\end{align*}
	we can choose $s$ large enough to absorb $Cm|v|^2$ by the dominant terms on the left-hand side of \eqref{Carlemaninw} and then we obtain the desired result for $L_q= \partial_t^2 -\text{div}(a(x)\nabla)+q$.

	\subsection{Carleman estimate with initial kinetic energy}


	\begin{corollary}[Carleman with initial kinetic energy]\label{Carliniitialwiths}
		Assuming the same conditions as stated in Theorem \ref{mainthm}, and with $\varphi$ given by \eqref{varphiweight}, if 
  $$v(0,x)=0, \quad \forall x\in \Omega,$$ 
  then there exists $s_0 > 0$ such that for all $s \geq s_0$, 
  	\begin{equation*}
		J(v):= s^{1/2}\int_{\Omega} e^{2s\varphi(0)} |\partial_t v(0)|^2 \,dx.
	\end{equation*}
 can be added to the left-hand side of \eqref{CarlBallR}.
	\end{corollary}

	Recalling the notation $w=e^{s\varphi }v$ and using the assumption $v(0,x)=0 \ \text{in}\ \Omega,$ we get $w(0,x)=0 \ \text{in}\ \Omega$. This last condition, combined with the cut-off function 
 $\chi \in C^\infty(\mathbb{R})$ such that $0\leq \chi\leq 1$ and
	\begin{align*}
		\chi(\tau)= 	\left\{ \begin{aligned} 
			1, \ & \quad \text{if} \quad \tau> 0,\\
			0, \ &\quad \text{if} \quad \tau \leq -T,
		\end{aligned}\right.
	\end{align*}
 allows us to write the following equivalences: 
	\begin{align}\label{wt0}
		\begin{split}
			s^{1/2}\int_{\Omega} e^{2s\varphi(0)}|\partial_t v(0)|^2\,dx&=	s^{1/2}\int_{\Omega} |\partial_t w(0)|^2\,dx\\ &= s^{1/2}\int_{\Omega}\left(\chi(0)|\partial_tw(0)|^2 - \chi(-T)|\partial_t w(-T)|^2 \right) \,dx \\ &= s^{1/2} \int_{\Omega} \int_{-T}^{0} \frac{d}{dt}\left(\chi(t)|\partial_t w(t)|^2 \right)dt\,dx\\
			&= s^{1/2} \int_{-T}^{0} \int_{\Omega} \left( \chi'(t)|\partial_t w(t)|^2 +2\chi(t)\partial_t w \partial_{t}^2 w \right)\,dxdt.
		\end{split}
	\end{align}
	Using the definition of $P_1 w$ given in \eqref{P1w}, we have
	$\partial_{t}^2 w = P_1w +a\Delta w -s^2(|\partial_t \varphi|^2 -a|\nabla \varphi|^2)w +2\nabla a\cdot \nabla w.$ Then, using Young's inequality,  \eqref{wt0} becomes
	\begin{align}\label{vtobynormw}
		\begin{split}
			&s^{1/2}\int_{\Omega} e^{2s\varphi(0)}|\partial_t v(0)|^2\,dx \leq\, Cs^{1/2} \int\!\!\!\int_{Q} |\partial_t w|^2 \,dxdt \\ &+ s^{1/2} \int_{-T}^{0} \int_{\Omega} 2\chi(t) \partial_t w \left(P_1w +a\Delta w -s^2(|\partial_t \varphi|^2 -a|\nabla \varphi|^2)w +2\nabla a\cdot \nabla w\right) \,dxdt\\
			&\leq C s \int\!\!\!\int_{Q} |\partial_t w|^2 \,dxdt + C\int\!\!\!\int_{Q} |P_1 w|^2 \,dxdt \\&\quad+ 2s^{1/2} \int_{-T}^{0} \int_{\Omega} \chi(t)\partial_t w (a\Delta w -s^2(|\partial_t \varphi|^2 -a|\nabla \varphi|^2)w +2\nabla a\cdot \nabla w)\,dxdt.
		\end{split}
	\end{align}
	The goal is to estimate the right-hand side of the previous inequality using Carleman estimates \eqref{Carlemaninw} in the variable $w$ . To achieve this, we need to perform integration by parts in the last integral. We can break down the sum of  this computation into three terms: $T_1, T_2$ and $T_3$.  Therefore, we can write the final result as 
	\begin{align}\label{wtP1w}
		2s^{1/2} \int_{-T}^{0}\int_{\Omega} \chi(t)\partial_t w (a\Delta w -s^2(|\partial_t \varphi|^2 -a|\nabla \varphi|^2)w +2\nabla a\cdot \nabla w)\,dx dt =	T_1 +T_2 +T_3.
	\end{align} 
	As is customary, we will begin by computing each term on a generic set $U$ without boundary conditions. We will then apply these computations to both sets $\Omega_1$ and $\Omega_2$, combining the common integrals supported on $\Gamma_*$. We can use the fact that $w(0,x)=0$ in $\Omega$ and $\chi(-T)=0$. By integrating by parts, using the transmission conditions \eqref{contrasw}, Young's inequality, and assumptions \eqref{infsup} - \eqref{reg}, it is easy to obtain the next estimate.
	\begin{align*}
		T_1 +T_2+ T_3 &\leq C  s\int\!\!\!\int_{Q}\left(|\partial_tw|^2 + a|\nabla w|^2 +s^2|w|^2\right) \,dxdt +  Cs^{3} \int\!\!\!\int_{\Sigma_*} |w|^2 \,d\sigma dt.
	\end{align*}

	Putting together the previous estimates for $T_1, T_2$ and $T_3$ in \eqref{wtP1w}, we can write the inequality \eqref{vtobynormw} as 
	
	\begin{align*}
		&s^{1/2}\int_{\Omega} e^{2s\varphi(0)}|\partial_t v(0)|^2\,dx\\ &\leq C  s\int\!\!\!\int_{Q}(|\partial_tw|^2 + a|\nabla w|^2 +s^2|w|^2) \,dxdt + 	s^3\int\!\!\!\int_{\Sigma_*} |w|^2\,d\sigma dt+ \int\!\!\!\int_{Q}|P_1w|^2\,dxdt,
	\end{align*}
	which using simultaneously \eqref{Carlemaninw} and \eqref{CarlBallR} give
	\begin{align}\label{usefulCarleman}
		\begin{aligned}
			s^{1/2}\int_{\Omega}& e^{2s\varphi(0)}|\partial_t v(0)|^2\,dx 	+s\int_{-T}^{T}\int_{\Omega}e^{2s\varphi} (|\partial_t v|^2 + a|\nabla v|^2 +s^2|v|^2)  \,dxdt \\
			&\leq C\int_{-T}^{T}\int_{\Omega}e^{2s\varphi} (L_qv)^2\,dxdt + Cs\int_{-T}^{T} \int_{\Gamma^+_{\mu^{x_1}}} e^{2s\varphi}\left|\partial_\nu v \right|^2 \,d\sigma dt\\ &\quad+ Cs \int\!\!\!\int_{\{(t,x) \in (-T,T)\times \Omega: \  \varphi<0\}}  e^{2s\varphi} (|\partial_t v|^2 +a|\nabla v|^2 + s^2|v|^2) \,dxdt \\
			&\quad +Cs\int_{-T}^{T}\int_{B_\varepsilon(x_1)}e^{2s\varphi} (|\partial_t v|^2 +a|\nabla v|^2 + s^2|v|^2) \,dxdt\\  &\quad+ Cs\int_{\Omega}e^{2s\varphi(\pm T)}(|\partial_t v(\pm T)|^2 +|\nabla v(\pm T)|^2 +s^2|v(\pm T)|^2) \,dx \\
			&\quad +Cs^3\int\!\!\!\int_{\{(t,x) \in (-T,T)\times\Gamma_*: \ \varphi<0\}} e^{2s\varphi}|v|^2\,d\sigma dt +Cs^2\int_{\Gamma_*} e^{2s\varphi(\pm T)}|v(\pm T)|^2 d\sigma.
		\end{aligned}
	\end{align}
	
	\section{Lipschitz Stability}\label{lipschitzsection}
	In this part of the paper, we present the proof of Theorem \ref{StabThm}. The proof strategy is to carefully follow the well-known approach for proving Lipschitz stability for the potential recovery, specifically in the scenario of smooth $a=a(x)$ wave speed. We will need to construct a convenient cut-off function that would enable us to apply estimate \eqref{usefulCarleman}. Furthermore, we will require the combination of two weights, each allowing us to remove localized integrals in $B_\varepsilon(x_1)$ and $B_\varepsilon(x_2)$.
	
	\subsection{Initiating the Bukgheim-Klibanov method}
	Let $u_p$ and $u_q$ be the corresponding solutions of \eqref{wavevardis} associated to potentials $p$ and~$q$. 
	Then
	\begin{align*}
		\begin{cases}
			(\partial^2_{t}  - \text{div} (a(x)\nabla) +q)(u_p-u_q) = (q-p)u_{p}, \quad &\text{in}  \quad  (0,T) \times \Omega ,  \\
			u_p-u_q = 0,   \quad &\text{on}  \quad   (0,T) \times \partial \Omega ,\\
			(u_p-u_q)(0) = 0, \quad \partial_t(u_p-u_q)(0) = 0, \quad &\text{in} \quad \Omega.
		\end{cases}
	\end{align*}
	
	Next we define $y=\partial_t (u_p-u_q)$ and we can notice that $ y$ satisfies 
	\begin{align}\label{wavevardisdifft}
		\begin{cases}
			(\partial^2_{t}  - \text{div} (a(x)\nabla) +q)y = (q-p)\partial_tu_{p}, \quad &\text{in}  \quad  (0,T) \times \Omega ,  \\
			y=0,   \quad &\text{on}  \quad   (0,T) \times \partial \Omega ,\\
			y(0)=0 \quad {\partial_t y}(0)=(q-p)u^0, \quad &\text{in} \quad \Omega.
		\end{cases}
	\end{align}
	
	Let us observe that the assumption $u_p \in H^1(0,T;L^\infty(\Omega))$ implies that $(q-p)\partial_tu_{p} \in L^1(0,T;L^2(\Omega))$ and $(q-p)u^0 \in L^2(\Omega)$. It is a classical result that under these assumptions and for $a=a(x)$ under the setting of assumptions \eqref{infsup}- \eqref{rega} and $q \in L^\infty (\Omega)$, equation \eqref{wavevardisdifft} admits a unique weak solution
	\begin{align}\label{regularityy}
		y \in C([0,T];H^1_0(\Omega)) \cap C^1([0,T];L^2(\Omega)),
	\end{align}
	where the energy $E_y(t)= \|\partial_t y (t)\|^2_{L^2(\Omega)} + a \|\nabla y (t)\|^2_{L^2(\Omega)}$ satisfies
	\begin{align}\label{energyEy}
		E_y(t) \leq C  \|(q-p)u^0\|^2_{L^2(\Omega)} +  C\|(q-p) \partial_t u_p\|^2_{L^1(0,T;L^2(\Omega))} .
	\end{align}
	Moreover, we can assert that $\partial_\nu y$ belongs to $L^2(0,T;L^2(\partial \Omega))$ and that there exists a positive constant $C$ depending in particular on $\|u_p\|_{ H^1(0,T;L^\infty(\Omega))}$ such that
	\begin{align*}
		\begin{aligned}
			\left\|\partial_\nu y\right\|^2_{L^2(0,T;L^2(\partial \Omega))} 
			& \leq C \|q-p\|^2_{L^2(\Omega)}.
		\end{aligned}
	\end{align*}

	To apply Corollary \ref{Carliniitialwiths} in this section, we need to extend the functions $y$ and $\partial_t u_{p}$ to $(-T,T)$ by taking their odd extensions in the $t$-variable. This odd extension enables us to maintain the regularity \eqref{regularityy} in $(-T,T)$. We will use the same notation for the extended functions. As a result, we obtain:
	\begin{align}\label{exty}
		\begin{cases}
			(\partial^2_{t}  - \text{div} (a(x)\nabla) +q)y = (q-p)\partial_tu_{p}, \quad &\text{in}  \quad  (-T,T) \times \Omega,  \\
			y=0,   \quad &\text{on}  \quad   (-T,T) \times \partial \Omega ,\\
			y(0)=0 \quad {\partial_t y}(0)=(q-p)u^0, \quad &\text{in} \quad \Omega.
		\end{cases}
	\end{align}
	
	The proof of Lipschitz stability can be obtained by estimating the initial velocity ${\partial_t y}(0)=(q-p)u^0$ in \eqref{exty} using the right-hand side of the same equation. To achieve this estimation, we can apply Carleman estimates \eqref{usefulCarleman}. However, it is important to ensure that the functions used in the Carleman estimate \eqref{usefulCarleman} meet the transmission conditions \eqref{transuj}. Moreover, these functions should allow us to eliminate additional observations that are supported on $(-T,T)\times B_\varepsilon(x_1)$, as well as terms evaluated at $\pm T$ and in the regions of $(-T,T)\times \Omega$ and $(-T,T)\times \Gamma_*$ where $\left\lbrace \varphi<0\right\rbrace$. Specifically, we will require:
	\begin{itemize}
		\item The combination of the weight $\varphi(t,x)=\mu^{x_1}(x)-\beta t^2$ with a second weight $\phi(t,x)= \mu^{x_2}(x) -\beta t^2$ defined in \eqref{phi} to eliminate the corresponding integrals supported on $(-T,T)\times B_\varepsilon(x_1)$ and $(-T,T)\times B_\varepsilon(x_2)$. 
		
		\item The construction of  a cut-off function $\theta^{\varphi}$ and $\theta^{\phi}$(depending either on $\varphi$ and $\phi$) to  get rid of pointwise terms evaluated at $\pm T$ and integrals supported on  $\left\lbrace \varphi<0\right\rbrace $ and $\left\lbrace \phi<0\right\rbrace$ respectively.
	\end{itemize}

	\subsection{Defining a convenient cut-off function}

	Consider
	\begin{equation*}
		\varphi(t,x)=\mu^{x_1}(x)-\beta t^2 \quad \text{and} \quad \phi(t,x)= \mu^{x_2}(x) - \beta t^2,
	\end{equation*}
	with $\mu^{x_1}$ and  $\mu^{x_2}$ previously defined. From \eqref{minkpositive} we can chose $\tilde{\varepsilon}< \min\{\delta_1,\delta_2\}$ with $\delta_1$ and $\delta_2$ satisfying 
	\begin{align*}
		\mu^{x_i}(x)\geq \delta_i >0 \quad \forall \,x \in (\Omega_1 \cup \Omega_2) \setminus B_\varepsilon(x_i), \qquad i=1,2.
	\end{align*}
	We can simply choose $\delta_1=\delta_2=\varepsilon^2 /[\text{diam}(\Omega)]^2$ and then take $\tilde{\varepsilon}<\varepsilon^2 /\left[ \text{diam}(\Omega)\right] ^2$ so that we have 
	\begin{align*}
		0<\tilde{\varepsilon}<\inf_{(\Omega_1 \cup \Omega_2) \setminus B_\varepsilon(x_k)} \mu^{x_k}(x), \quad k=1,2.
	\end{align*}
	
	Recalling the notation $\varphi_j(t,x)= \mu^{x_1}\mathbf{1}_{\Omega_j}(x) -\beta t^2$  for $j=1,2$, we define  $\theta^\varphi \in C^\infty((-T,T);\Omega)$ as follows.  First, in $(-T,T)\times \Omega_1$ we define

	\begin{align}
		\theta_1^{\varphi}(t,x)= 	\left\{ \begin{aligned} \label{theta1varphi}
			0, \ & \quad \text{if} \quad \varphi_1(t,x)< M_1,\\
			1, \ &\quad \text{if} \quad \varphi_1 (t,x)\geq M_1 + \tilde{\varepsilon}
		\end{aligned}\right.
	\end{align}
	and in $(-T,T)\times \Omega_2$
	
	\begin{align}
		\theta_2^{\varphi}(t,x)= 	\left\{ \begin{aligned} \label{theta2varphi}
			0, \ & \quad \text{if} \quad \varphi_2(t,x)< M_2,\\
			1, \ &\quad \text{if} \quad \varphi_2 (t,x) \geq M_2 + \gamma,
		\end{aligned}\right.
	\end{align}
	for constants $M_1$, $M_2$  and $\gamma$ given in  \eqref{constantsM1M2x1gamma}. The key observation is that the variation from zero to one must be done in $(-T,T)\times \Omega_1$ and $(-T,T)\times \Omega_2$  satisfying particularly the trace condition  
	\begin{align}\label{transthetavarphi}
		\frac{\partial \theta_1^{\varphi}}{\partial \nu_1} = \frac{\partial \theta_2^{\varphi}}{\partial \nu_2} = 0, \quad &\text{on} \quad (0,T) \times \Gamma_*. 
	\end{align}
	
	Note that for the second weight, $\phi(t,x)= \mu^{x_2}(x) - \beta t^2$, we  define $\theta ^{\phi} \in C^\infty((-T,T)\times \Omega)$ in a strictly analogous fashion of \eqref{theta1varphi} and \eqref{theta2varphi}.  
	
\begin{figure}[ht]\label{tracestheta}
\centering
\includegraphics[width=1\textwidth]{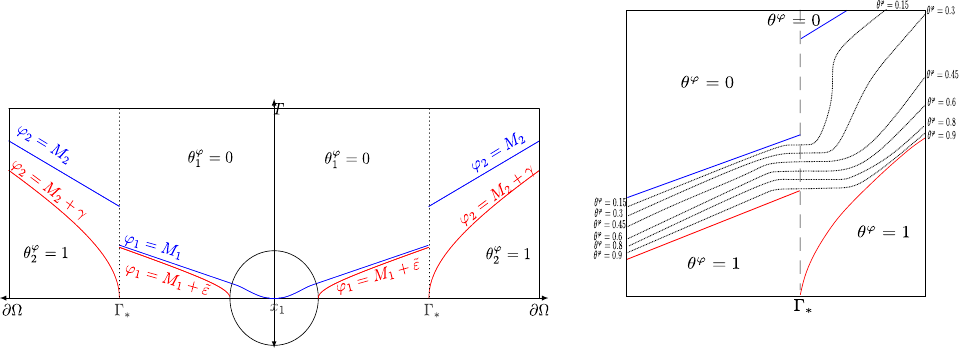}
\caption{Context of the application of $\theta^{\varphi}$ and the behaviour of $\theta^{\varphi}$ in a neighborhood of $\Gamma_*$ satisfying \eqref{transthetavarphi}.}
\end{figure}

Figure \ref{tracestheta}(left) represents the context of the application of $\theta^{\varphi}$.  The variation of $\theta_1^{\varphi}$ in $(-T,T)\times \Omega_1$ needs to be done between the upper blue line $(\varphi_1= M_1)$ and the lower red curved line $(\varphi_1= M_1+\tilde{\varepsilon})$. Similarly, in $(-T,T)\times \Omega_2$, the variation of $\theta_2^{\varphi}$ must be done between the upper blue line $(\varphi_2=M_2)$ and the lower red curve line $(\varphi_2= M_2 +\gamma)$.
	
In the $1D$ case, condition \eqref{transthetavarphi} can be achieved by allowing just a “perpendicular variation" at the points of discontinuity. See in particular  \cite[Section 3.1]{baudouin:hal-04361363} for the situation of two connected branches in the case of waves on networks. In the setting of Figure \ref{tracestheta}(right), it means that the cut-off function's level sets going from $0$ above the blue line toward $1$ below the red line are parallel when crossing $\Gamma_*$ perpendicularly.
	
	
	\subsection{Applying Carleman estimates}
	
	Let us take $\beta$ satisfying 
	\begin{equation}\label{betaapplication}
		\frac{L^2}{T^2}<\beta<{\beta}_0,
	\end{equation}
	which can be done since the assumption \eqref{timecond} on $T$. Then, using the function $\theta^{\varphi}$ defined  in \eqref{theta1varphi} and \eqref{theta2varphi} and the analogous function $\theta^{\phi}$ 
 let us set
	\begin{align*}
		v^1 = \theta^{\varphi}y \quad \text{and} \quad v^2 = \theta^{\phi}y,
	\end{align*}
	for $y$ solving \eqref{exty}. Note that as soon as condition \eqref{transthetavarphi} is satisfied for $\theta^\varphi$ and $ \theta^{\phi}$ respectively and in view of the condition  for $\beta$ in \eqref{betaapplication}, then $v^1$ and $v^2$ satisfies the next properties.
	
	\begin{itemize}
		\item $v^1$ and $v^2$ satisfies  the same transmission conditions as \eqref{transuj}:
		\begin{align*}
			\begin{aligned}
				v_1^1 = v_2^1, \quad \text{and} \quad  a_1\frac{\partial v_1^1}{\partial \nu_1} + a_2\frac{\partial v_2^1}{\partial \nu_2} = 0, \quad &\text{on} \quad (0,T) \times \Gamma_*,\\
				v_1^2 = v_2^2, \quad \text{and} \quad  a_1\frac{\partial v_1^2}{\partial \nu_1} + a_2\frac{\partial v_2^2}{\partial \nu_2} = 0, \quad &\text{on} \quad (0,T) \times \Gamma_*.
			\end{aligned}
		\end{align*}
		\item $v^1(\pm T)=v^2(\pm T)=0 $ and  ${\partial_t v^1}(\pm T)={\partial_t v^2}(\pm T)=0, \ \text{in} \ \Omega$.
		
		\item $v^1$ and its derivatives vanishes both in the region $\{(t,x) \in (-T,T) \times \Omega : \ \varphi < 0\}$  and $\{(t,x) \in (-T,T) \times \Gamma_* : \ \varphi < 0\}$. The same applies to $v^2$ in the corresponding regions $\{(t,x) \in (-T,T) \times \Omega : \ \phi < 0\}$ and  $\{(t,x) \in (-T,T) \times \Gamma_* : \ \phi < 0\}$.
		
		\item  The equations
		\begin{align}\label{wavevardisv1}
			\begin{cases}
				(\partial^2_{t}  - \text{div} (a(x)\nabla) +q)v^1 = \theta^\varphi(q-p)\partial_tu_{p} + [L_q,\theta^\varphi]y, \quad &\text{in}  \quad  (-T,T) \times \Omega ,  \\
				v^1=0,   \quad &\text{on}  \quad   (-T,T) \times \partial \Omega ,\\
				v^1(0)=0, \quad {\partial_t v^1}(0)=\theta^{\varphi}(0)(q-p)u_p(0), \quad &\text{in} \quad \Omega,\\
			\end{cases}
		\end{align}
		and 
		\begin{align}\label{wavevardisv2}
			\begin{cases}
				(\partial^2_{t}  - \text{div} (a(x)\nabla) +q)v^2 = \theta^\phi(q-p)\partial_tu_{p} + [L_q,\theta^\phi]y, \quad &\text{in}  \quad  (-T,T) \times \Omega ,  \\
				v^2=0,   \quad &\text{on}  \quad   (-T,T) \times \partial \Omega ,\\
				v^2(0)=0, \quad {\partial_t v^2}(0)=\theta^\phi(0)(q-p)u_p(0), \quad &\text{in} \quad \Omega,
			\end{cases}
		\end{align}
		are satisfied. Here $[A,B]y=A(By)- B(Ay)$ stands for the commutator of two operators $A$ and $B$. Additionally $a=a(x)$ also satisfies the multiplier condition \eqref{rega}  and  $q \in L^\infty_{\leq m}(\Omega)$.
	\end{itemize}
	
	Now, we are in the position of applying the inequality  \eqref{usefulCarleman} to both $v^1$ and $v^2$ and to sum up the resulting inequalities. Consequently, we obtain, for all $s\geq s_0$,
	
	\begin{align}\label{usefulCarlemanv1v2}
		\begin{aligned}
			&s^{1/2}\int_{\Omega} e^{2s\varphi(0)} |\partial_t v^1(0)|^2  \,dx \, + 
			s\int_{-T}^{T}\int_{\Omega}e^{2s\varphi} (|\partial_t v^1|^2 +a|\nabla v^1|^2 + s^2|v^1|^2) \,dxdt \\
			& +s^{1/2}\int_{\Omega} e^{2s\phi(0)} |\partial_t v^2(0)|^2 \,dx +  s\int_{-T}^{T}\int_{\Omega}e^{2s\phi} (|\partial_t v^2|^2 +a|\nabla v^2|^2 + s^2|v^2|^2) \,dxdt \\ 
			&\leq C\int_{-T}^{T}\int_{\Omega}e^{2s\varphi} (L_qv^1)^2\,dxdt + C\int_{-T}^{T}\int_{\Omega}e^{2s\phi} (L_qv^2)^2\,dxdt\\ &\quad+ Cs\int_{-T}^{T} \int_{\Gamma^+_{\mu^{x_1}}} e^{2s\varphi}\left|\partial_\nu v^1 \right|^2 \,d\sigma dt+ Cs\int_{-T}^{T}\int_{\Gamma^+_{\mu^{x_2}}} e^{2s\phi}\left|\partial_\nu v^2 \right|^2 \,d\sigma dt\\
			&\quad+Cs\int_{-T}^{T}\int_{B_\varepsilon(x_1)}e^{2s\varphi} (|\partial_t v^1|^2 +a|\nabla v^1|^2 + s^2|v^1|^2) \,dxdt\\
			&\quad+Cs\int_{-T}^{T}\int_{B_\varepsilon(x_2)}e^{2s\phi} (|\partial_t v^2|^2 +a|\nabla v^2|^2 + s^2|v^2|^2) \,dxdt.
		\end{aligned}
	\end{align}

	\subsection{Eliminating ``blinds spots" integrals}
	In this step, we aim to eliminate the last two integrals (on $B_\varepsilon(x_1)$ and $B_\varepsilon(x_2)$) of the right-hand side of estimate \eqref{usefulCarlemanv1v2}. We will do this in two parts. Firstly, we will compare $v^1$ and its derivatives with $v^2$ in the region $(-T,T) \times B_\varepsilon(x_1)$. Simultaneously, we will also compare $v^2$ with $v^1$ in the region $(-T,T) \times B_\varepsilon(x_2)$. For this comparison to be possible, the value of $\varepsilon$ needs to be small enough. 
	
	In the second part, we will use the integrals on the left-hand side of equation \eqref{usefulCarlemanv1v2} to absorb the integrals on $B_\varepsilon(x_1)$ and $B_\varepsilon(x_2)$. This absorption will be possible by choosing an appropriate value of $s$ that is large enough.
	
	Note that using \eqref{minkbmo} we have the estimate
	\begin{align}\label{varphiinBx1}
		\varphi(t,x) = \mu^{x_1}(x) - \beta t^2 &\leq \left[ p_{\Omega_1^{x_1}}(x-x_1)\right]^2 - \beta t^2 +M_1\nonumber \\ &\leq \frac{\varepsilon^2}{\left[ \text{dist}(x_1,\Gamma_*)\right]^2 } - \beta t^2  +M_1\quad \text{in } (-T,T)\times B_\varepsilon(x_1).
	\end{align}
	In addition, if we assume  
	\begin{align}\label{epsilond}
		\varepsilon < d:=\frac{|x_1 - x_2|}{2}
	\end{align}
	 we also have
	
	\begin{align}\label{phibelowBx1}
		\phi(t,x) &= \left[ p_{\Omega_1^{x_2}}(x-x_2) \right]^2 - \beta t^2 +M_1= \frac{|x-x_2|^2}{|y(x)-x_2|^2} - \beta t^2 +M_1 \nonumber
		\\ & \geq \frac{d^2}{\left[ \max_{y \in \Gamma_*}\text{dist}(y, x_2)\right]^2 }  - \beta t^2 +M_1 \quad \text{in } (-T,T)\times B_\varepsilon(x_1). 
	\end{align}
	
	Notice that the right-hand side of \eqref{varphiinBx1} depends on $t \in (-T,T)$. We are particularly interested in the points of $t\in (-T,T)$ where 
	\begin{align}\label{varphibelowBx1sign}
		\frac{\varepsilon^2}{\left[ \text{dist}(x_1,\Gamma_*)\right]^2 } - \beta t^2  +M_1 \geq M_1,
	\end{align}
	because in the other case, the function  $v^1$ vanishes according to the definition of $\theta^\varphi$.  
	Using \eqref{varphibelowBx1sign} in \eqref{phibelowBx1} we obtain
	
	\begin{align*}
		\phi(t,x) \geq \frac{d^2}{\left[ \max_{y \in \Gamma_*}\text{dist}(y, x_2)\right]^2 }  - \frac{\varepsilon^2}{\left[ \text{dist}(x_1,\Gamma_*)\right]^2 }+M_1 \quad \text{in } (-T,T)\times B_\varepsilon(x_1),
	\end{align*}
	
	which for
	
	\begin{align}\label{epsilonsmall1}
		\varepsilon< \frac{d}{\max_{y \in \Gamma_*}\text{dist}(y, x_2) } \frac{\text{diam} (\Omega)\, \text{dist}(x_1,\Gamma_*) }{\sqrt{\left[ \text{dist}(x_1,\Gamma_*)\right]^2 + \left[ \text{diam}(\Omega)\right]^2 }}:=d_1,
	\end{align}
	implies 
	\begin{align} \label{phigreatertildepsi}
		\phi(t,x) \geq M_1 +  \tilde{\varepsilon} \quad \text{in } (-T,T)\times B_\varepsilon(x_1)
	\end{align}
	and consequently
	\begin{align}\label{v1inball}
		|v^1|^2 = |\theta^\varphi y|^2 \leq |y|^2 \leq |\theta^\phi y|^2 = |v^2|^2 \quad   \text{in} \ (-T,T)\times B_\varepsilon(x_1).
	\end{align}
	
	It is direct to check that the condition  \eqref{phigreatertildepsi}  implies that $\nabla \theta^\phi$ and $\partial_t \theta^\phi$ neglects in $(-T,T)\times B_\varepsilon(x_1) $. Then we obtain
	\begin{align*}
		|\nabla v^2|^2 =|\nabla \theta^{\phi} y + \theta^\phi \nabla y|^2=|\nabla y|^2   \quad \text{in } (-T,T)\times B_\varepsilon(x_1)
	\end{align*}
	and
	\begin{align*}
		|\partial_t v^2|^2 =|\partial_t \theta^{\phi} y + \theta^\phi \partial_t y|^2=|\partial_t y|^2   \quad \text{in } (-T,T)\times B_\varepsilon(x_1).
	\end{align*}
	These two last conditions can be used together with \eqref{v1inball} to get 
	\begin{align}\label{gradviball}
		|\nabla v^1|^2 &=  |\nabla \theta^\varphi y + \theta^\varphi \nabla y|^2 \nonumber\\ &\leq C|y|^2 + C|\nabla y|^2\leq C|v^2|^2 + C|\nabla v^2|^2   \quad \text{in } (-T,T)\times B_\varepsilon(x_1)
	\end{align}
	and 
	\begin{align}\label{v1tinball}
		|\partial_t v^1|^2 \leq C|v^2|^2 + C|\partial_t v^2|^2  \quad \text{in } (-T,T)\times B_\varepsilon(x_1).
	\end{align}

	Combining \eqref{v1inball}, \eqref{gradviball} and \eqref{v1tinball}, we observe that following estimate holds
	\begin{align}\label{v1byv2in0}
		\begin{aligned}
			s\int_{-T}^{T}\int_{B_\varepsilon(x_1)}e^{2s\varphi} (|\partial_t v^1|^2 &+a|\nabla v^1|^2 + s^2|v^1|^2) \,dx dt \\
			\leq &Cs\int_{-T}^{T}\int_{B_\varepsilon(x_1)}e^{2s\varphi} (|\partial_t v^2|^2 +a|\nabla v^2|^2 + s^2|v^2|^2) \,dx dt.
		\end{aligned}
	\end{align}

	Performing  analogous steps to arrive to \eqref{v1byv2in0}, it is straightforward to obtain an equivalent inequality to \eqref{v1byv2in0} in  $(-T,T) \times B_\varepsilon (x_2)$; that is, if 
	
	\begin{align}\label{epsilonsmall2}
		\varepsilon< \frac{d}{\max_{y \in \Gamma_*}\text{dist}(y, x_1) } \frac{\text{diam} (\Omega)\, \text{dist}(x_2,\Gamma_*) }{\sqrt{\left[ \text{dist}(x_2,\Gamma_*)\right]^2 + \left[ \text{diam}(\Omega)\right]^2 }}:=d_2,
	\end{align}
	then
	\begin{align}\label{v2byv1inx2}
		\begin{aligned}
			s\int_{-T}^{T}\int_{B_\varepsilon(x_2)}e^{2s\phi} (|\partial_t v^2|^2 +&a|\nabla v^2|^2 + s^2|v^2|^2) \,dx dt \\
			\leq &Cs\int_{-T}^{T}\int_{B_\varepsilon(x_2)}e^{2s\phi} (|\partial_t v^1|^2 +a|\nabla v^1|^2 + s^2|v^1|^2) \,dx dt.
		\end{aligned}
	\end{align}

	Let us fix $\varepsilon>0$ such that \eqref{epsilond},\eqref{epsilonsmall1} and \eqref{epsilonsmall2} hold, i.e.,
	
	\begin{align}\label{epsilonmincondition}
		\varepsilon<\min\left\lbrace d, d_1,d_2 \right\rbrace
	\end{align}
	and note that in particular that the condition \eqref{epsilonmincondition}  also implies
	\begin{align}\label{differenceofweightsx1}
		\phi (t,x) - \varphi(t,x) = \mu^{x_2}(x)-\mu^{x_1}(x) \geq \tilde{\varepsilon} > 0 \quad \text{in}  \quad (-T,T)\times B_\varepsilon(x_1)
	\end{align}
	and
	\begin{align}\label{differenceofweightx2}
		\varphi (t,x) - \phi(t,x) = \mu^{x_1}(x)-\mu^{x_2}(x) \geq \tilde{\varepsilon} > 0 \quad \text{in}  \quad (-T,T)\times B_\varepsilon(x_2).
	\end{align}
	Consequently, we observe that 
	\begin{align}\label{finalstepballex1}
		Cs e^{2s\varphi} \leq \frac{s}{2} e^{2s\phi},  \quad \text{in}  \quad  (-T,T)\times B_\varepsilon(x_1)
	\end{align}
	and
	\begin{align}\label{finalstepballex2}
		Cs e^{2s\phi} \leq \frac{s}{2} e^{2s\varphi},  \quad \text{in}  \quad  (-T,T)\times B_\varepsilon(x_2),
	\end{align}
	for all $s\geq s_0$. 
	
	Plugging \eqref{finalstepballex1} and \eqref{finalstepballex2} in \eqref{v1byv2in0} and \eqref{v2byv1inx2} respectively we can finally take  $s$ sufficiently large to obtain (thanks to inequality \eqref{usefulCarlemanv1v2})
	\begin{align}\label{usefulCarlemanv1v2fin}
		\begin{aligned}
			s^{1/2}\int_{\Omega} &e^{2s\varphi(0)} |\partial_t v^1(0)|^2 \,dx + s^{1/2}\int_{\Omega} e^{2s\phi(0)} |\partial_t v^2(0)|^2 \,dx\\
			&\leq C\int_{-T}^{T}\int_{\Omega}e^{2s\varphi} (L_qv^1)^2 \,dxdt + C\int_{-T}^{T}\int_{\Omega}e^{2s\phi} (L_qv^2)^2 \,dxdt\\ 
			&+ Cs\int_{-T}^{T} \int_{\Gamma^+_{\mu^{x_1}}} e^{2s\varphi}\left|\partial_\nu v^1 \right|^2 \,d\sigma dt+ Cs\int_{-T}^{T}\int_{\Gamma^+_{\mu^{x_2}}} e^{2s\phi}\left|\partial_\nu v^2 \right|^2 \,d\sigma dt.
		\end{aligned}
	\end{align}
	
	\subsection{Applying  Bukhgeim-Klibanov method }
	
	Once we obtain \eqref{usefulCarlemanv1v2fin}, the strategy for proving Lipschitz stability in the inverse problem is rather classic. We will need to pay specific attention to carefully combine the information from the left-hand side of \eqref{usefulCarlemanv1v2fin} for getting an estimate of $\|p-q\|_{L^2(\Omega)}$. 
	
	Taking into account the definition of $\theta^\varphi$ and $\theta^\phi$, it is clear that 
	\begin{align*}
		\theta^{\varphi}(0,x)=1,  \quad \text{in}  \quad  \Omega\setminus B_\varepsilon(x_1),\\
		\theta^{\phi}(0,x)=1,  \quad \text{in}  \quad  \Omega\setminus B_\varepsilon(x_2).
	\end{align*}
	Using the assumption $|u^0(x)|\geq \delta >0$ a.e. in $\Omega$, the left hand side of \eqref{usefulCarlemanv1v2fin}  is minimized as follows 
	
	\begin{align}\label{firststepstabilityleft}
		\begin{aligned}
			s^{1/2}\delta^2\int_{\Omega\setminus B_\varepsilon(x_1)} e^{2s\varphi(0)} |q-p|^2   \,dx+ s^{1/2}\delta^2\int_{\Omega\setminus B_\varepsilon(x_2)} e^{2s\phi(0)} |q-p|^2  \,dx \\ \leq C 	s^{1/2}\int_{\Omega} e^{2s\varphi(0)} |\partial_t v^1(0)|^2 \,dx + s^{1/2}\int_{\Omega} e^{2s\phi(0)} |\partial_t v^2(0)|^2 \,dx.
		\end{aligned}
	\end{align}
	
	We now deal with the integrals of the right hand side of \eqref{usefulCarlemanv1v2fin} involving $|L_q v^1|^2$ and $|L_q v^2|^2$. From \eqref{wavevardisv1} and \eqref{wavevardisv2} these terms satisfy
	
	\begin{align}
		|L_q v^1|^2 \leq 2|\theta^\varphi (q-p)\partial_t u_p|^2 +2\left| [L_q,\theta^\varphi]y\right| ^2, \label{Lqv1sinint}\\
		|L_q v^2|^2 \leq 2|\theta^\phi (q-p)\partial_t u_p|^2 +2\left| [L_q,\theta^\phi]y\right| ^2. \label{Lqv2sinint}
	\end{align}
	
	It is immediate to note that the commutator $[L_q,\theta^\varphi]$ is confined in regions depending upon $\varphi$.  First, in $(-T,T)\times\Omega_1$ and $(-T,T)\times\Omega_2$ we see that $[L_q,\theta^\varphi]$ is confined at the points  $\left\lbrace M_1< \varphi_1< M_1 + \tilde{\varepsilon}  \right\rbrace $ and $\left\lbrace M_2< \varphi_2< M_2 + \gamma \right\rbrace $ respectively.   Therefore, substituting \eqref{Lqv1sinint} into the first integral of the right hand side of \eqref{usefulCarlemanv1v2fin} and since $\varphi(t,x)\leq \varphi(0,x)$ in $ (-T,T)\times \Omega$ we obtain
	
	\begin{align}\label{Lqv1}
		\begin{aligned}
			\int_{-T}^{T}\int_{\Omega}e^{2s\varphi} (L_qv^1)^2\,dxdt \leq &C \int_{-T}^{T}\int_{\Omega}e^{2s\varphi(0)} |q-p|^2|\partial_t u_p|^2 \,dx\,dt \\ &+ Ce^{2s(M_1+ \tilde{\varepsilon})}\int_{-T}^{T}\int_{\Omega_1} (|y|^2 +a_1|\nabla y|^2 +|\partial_ty|^2)\,dxdt
			\\ &+ Ce^{2s({M_2 +\gamma})}\int_{-T}^{T}\int_{\Omega_2} (|y|^2 +a_2|\nabla y|^2 +|\partial_ty|^2)\,dxdt.
		\end{aligned}
	\end{align}
	
	Let us return to the interval $(0,T)$ in \eqref{Lqv1}.To estimate $|y|^2$, we can use the Poincaré inequality. Then, we can use the energy estimates in equation \eqref{energyEy} for $y \in H_0^1(\Omega)$ in the last two integrals of the right-hand side of equation \eqref{Lqv1}. This allows us to conclude that there exists a positive constant (which  change from line to line) depending on $\|u_p\|_{H^1(0,T;L^\infty(\Omega))}$, $T$, and $m$ such that 
	
	\begin{align}\label{Lqv10T}
		\begin{aligned}
			\int_{-T}^{T}\int_{\Omega}e^{2s\varphi} (L_qv^1)^2\,dx dt\leq &C \int_{\Omega} e^{2s\varphi(0)} |q-p|^2\,dx + Ce^{2s(M_1+ \tilde{\varepsilon})}\int_{\Omega_1} |q-p|^2\,dx  \\ &+ Ce^{2s(M_2+ \gamma)}\int_{\Omega_2} |q-p|^2\,dx.
		\end{aligned}
	\end{align}
	
	Moreover, according to the estimates
	\begin{align*}
		\varphi_1(0,x)= \frac{|x-x_1|^2}{|y(x)-x_1|^2} +M_1  \geq \frac{\varepsilon^2}{\left[\text{diam}(\Omega) \right]^2  } +M_1  > \tilde{\varepsilon} +M_1 \quad \text{in} \quad \Omega_1 \setminus B_\varepsilon(x_1)
	\end{align*}
	and
	\begin{align*}
		\varphi_2(0,x)= \gamma \frac{|x-x_1|^2}{|y(x)-x_1|^2} +M_2 \geq \gamma +M_2  \quad \text{in} \quad \Omega_2,
	\end{align*}
	we  deduce from \eqref{Lqv10T} that 
	
	\begin{align}\label{Lqv1Final}
		\int_{-T}^{T}\int_{\Omega}e^{2s\varphi} (L_qv^1)^2\,dxdt &\leq C \int_{\Omega\setminus B_\varepsilon(x_1)} e^{2s\varphi(0)} |q-p|^2\,dx \nonumber \\ &\quad+  C\int_{B_\varepsilon(x_1)}\left( e^{2s(M_1+ \tilde{\varepsilon})} + e^{2s\varphi(0)}\right)  |q-p|^2\,dx.
	\end{align}
	
	Similar arguments can be applied to the second integral in the right hand side of \eqref{usefulCarlemanv1v2fin}. More precisely, using \eqref{Lqv2sinint} together with the confinedness  property of $[L_q,\theta^\phi]$ in the corresponding regions and since
	\begin{align*}
		\phi_1(0,x)= \frac{|x-x_2|^2}{|y(x)-x_2|^2} +M_1  \geq \frac{\varepsilon^2}{\left[\text{diam}(\Omega) \right]^2  } +M_1  > \tilde{\varepsilon} +M_1 \quad \text{in} \quad \Omega_1 \setminus B_\varepsilon(x_2)
	\end{align*}
	and
	\begin{align*}
		\phi_2(0,x)= \gamma \frac{|x-x_2|^2}{|y(x)-x_2|^2} +M_2 \geq \gamma +M_2  \quad \text{in} \quad \Omega_2,
	\end{align*}
	we also conclude that  there exists a positive constant that may change from line to line and depending in particular of $\|u_p\|_{H^1(0,T;L^\infty(\Omega))}, T$ and $m$ such that 
	\begin{multline*}
	\int_{-T}^{T}\int_{\Omega}e^{2s\phi} (L_qv^2)^2\,dxdt 
 \leq C \int_{\Omega\setminus B_\varepsilon(x_2)} e^{2s\phi(0)} |q-p|^2\,dx \nonumber \\ 
  \quad+  C\int_{B_\varepsilon(x_2)}\left( e^{2s(M_1+ \tilde{\varepsilon})} + e^{2s\phi(0)}\right)  |q-p|^2\,dx.
	\end{multline*}
	
	After returning to the interval $(0,T)$ and using the boundary conditions of $v^1$ and $v^2$ in \eqref{wavevardisv1} and \eqref{wavevardisv2}, the boundary integrals in \eqref{usefulCarlemanv1v2fin} can be expressed as follows in terms of the variable $y=\partial_t(u_p-u_q)$:
	
	\begin{multline*}
	    s\int_{-T}^{T} \int_{\Gamma_{\mu^{x_1}}} e^{2s\varphi}\left|\partial_\nu v^1 \right|^2 \,d\sigma dt+ s\int_{-T}^{T}\int_{\Gamma_{\mu^{x_2}}} e^{2s\phi}\left|\partial_\nu v^2 \right|^2 \,d\sigma dt \\ 
     \leq 
	Cs\int_{0}^{T} \hspace{-0.2cm}\int_{\Gamma^+_{\mu^{x_1}}} e^{2s\varphi}\left|\partial_\nu \partial_t (u_p -u_q)\right|^2 \,d\sigma dt
   + Cs\int_{0}^{T}\hspace{-0.2cm}\int_{\Gamma^+_{\mu^{x_2}}} e^{2s\phi}\left|\partial_{\nu} \partial_t (u_p-u_q) \right|^2 \,d\sigma dt.
	\end{multline*}
	
	Summarizing, we have that using \eqref{firststepstabilityleft}, \eqref{Lqv1Final}, and these last two estimates, the inequality \eqref{usefulCarlemanv1v2fin} becomes
	
	\begin{align}\label{almostStability}
		\begin{aligned}
			&s^{1/2}\int_{\Omega\setminus B_\varepsilon(x_1)}e^{2s\varphi(0)} |q-p|^2   dx+ s^{1/2}\int_{\Omega\setminus B_\varepsilon(x_2)} e^{2s\phi(0)} |q-p|^2 dx \\ 
   &\leq C\int_{\Omega\setminus B_\varepsilon(x_1)} e^{2s\varphi(0)} |q-p|^2 dx 
   + C \int_{\Omega\setminus B_\varepsilon(x_2)} e^{2s\phi(0)} |q-p|^2 dx \\ & +C\int_{B_\varepsilon(x_1)}\hspace{-0.2cm}\left( e^{2s(M_1+ \tilde{\varepsilon})} + e^{2s\varphi(0)}\right)  |q-p|^2 dx +  C\int_{B_\varepsilon(x_2)}\hspace{-0.2cm}\left( e^{2s(M_1+ \tilde{\varepsilon})} + e^{2s\phi(0)}\right)  |q-p|^2 dx\\
			& + Cs\int_{0}^{T} \int_{\Gamma^+_{\mu^{x_1}}} e^{2s\varphi}\left|\partial_\nu \partial_t (u_p -u_q)\right|^2 d\sigma dt
   + Cs\int_{0}^{T}\int_{\Gamma^+_{\mu^{x_2}}} e^{2s\phi}\left|\partial_{\nu} \partial_t (u_p-u_q) \right|^2 d\sigma dt.
		\end{aligned}
	\end{align}

	Next, we notice that taking $s$ sufficiently large, the first two integrals in the right-hand side of \eqref{almostStability} can be absorbed by the two integrals in the left-hand side of  \eqref{almostStability}. Furthermore, we note that
	\begin{align*}
		\phi(0,x)&= \frac{|x-x_2|^2}{|y(x)-x_2|^2} +M_1  \geq \frac{\varepsilon^2}{\left[\text{diam}(\Omega) \right]^2  } +M_1  > M_1+ \tilde{\varepsilon}  \quad \text{in} \quad B_\varepsilon(x_1)
	\end{align*}
	and from \eqref{differenceofweightsx1} we also see that
	
	\begin{equation*}
		\varphi(0,x) \leq \phi(0,x) \quad \text{in} \quad B_\varepsilon(x_1).
	\end{equation*}
	Analogously, we can assert that
	\begin{align*}
		\varphi(0,x)&= \frac{|x-x_1|^2}{|y(x)-x_1|^2} +M_1  \geq \frac{\varepsilon^2}{\left[\text{diam}(\Omega) \right]^2  } +M_1  > M_1+\tilde{\varepsilon}  \quad \text{in} \quad B_\varepsilon(x_2)
	\end{align*}
	and from \eqref{differenceofweightx2} we obtain
	\begin{equation*}
		\phi(0,x) \leq \varphi(0,x) \quad \text{in} \quad B_\varepsilon(x_2).
	\end{equation*}
	According to the latest estimates and since $\varepsilon$ satisfies \eqref{epsilonmincondition}, the third and fourth integrals on the right-hand side of \eqref{almostStability} satisfy
	$$
	    \int_{B_\varepsilon(x_1)}\left( e^{2s(M_1+ \tilde{\varepsilon})} + e^{2s\varphi(0)}\right)  |q-p|^2\,dx  \leq C\int_{\Omega \setminus B_\varepsilon(x_2)}e^{2s\phi(0)} |q-p|^2\,dx
	$$
	and
	$$
	    \int_{B_\varepsilon(x_2)}\left( e^{2s(M_1+ \tilde{\varepsilon})} + e^{2s\phi(0)}\right)  |q-p|^2\,dx \leq C\int_{\Omega \setminus B_\varepsilon(x_1)}e^{2s\varphi(0)} |q-p|^2\,dx.
	$$
	Taking $s$ large enough in the left hand side of  \eqref{almostStability} to absorb the integrals in the right hand side of these last two estimates we finally obtain
	\begin{align}\label{almstab2}
		\begin{aligned}
			&s^{1/2} \int_{\Omega\setminus B_\varepsilon(x_1)}e^{2s\varphi(0)} |q-p|^2 \,dx + s^{1/2} \int_{\Omega\setminus B_\varepsilon(x_2)}e^{2s\phi(0)} |q-p|^2 \,dx \\ &\leq
			Cs\int_{0}^{T} \int_{\Gamma^+_{\mu^{x_1}}} e^{2s\varphi}\left|\partial_\nu \partial_t (u_p -u_q)\right|^2 \,d\sigma dt+ Cs\int_{0}^{T}\int_{\Gamma^+_{\mu^{x_2}}} e^{2s\phi}\left|\partial_{\nu} \partial_t (u_p-u_q) \right|^2 \,d\sigma dt.    
		\end{aligned}
	\end{align}

	It follows immediately that for fixed $s$, we can estimate the exponential on both sides of \eqref{almstab2}. Indeed, using that we can cover $\Omega$ by $\Omega\setminus B_{\varepsilon}(x_1)$ and $\Omega\setminus B_{\varepsilon}(x_2)$, we finally obtain, for $\Gamma_0=\Gamma^+_{\mu^{x_1}}\cup \Gamma^+_{\mu^{x_2}}$, the inequality
	\begin{align*}
		\int_{\Omega}|q-p|^2 \,dx \leq C \int_{0}^{T}\int_{\Gamma_0} \left|\partial_{\nu} \partial_t (u_p-u_q) \right|^2 \,d\sigma dt.
	\end{align*}
	This ends the proof of Theorem \ref{StabThm} detailing the Lipschitz stability of the inverse problem considered in this article. 
	
	\ack
   A. Imba was partially supported by ANID BECAS/DOCTORADO NACIONAL 21221608. A. Mercado-Saucedo was
 partially supported by Fondecyt 1211292 and	ANID Millennium Science Initiative Program, Code NCN19-161,
 and A. Osses was partially supported by CMM FB210005 Basal-ANID, ANID-Fondecyt 1240200, 1231404, Fondef IT2310095
and FONDAP/15110009 grants.

	\section*{References}
	\bibliographystyle{plain}

\end{document}